%% file: ClassifiactionOfChannelCMCSurfaces.tex
\documentclass[a4paper,12 pt,final]{amsart} 

\input{Preamble/Packages}
\input{Preamble/Environments}
\input{Preamble/Commands}

\begin{document}
\title{Rotational cmc surfaces in terms of Jacobi elliptic 
functions}
\author{Denis Polly}

\date{\today}

\begin{abstract}
	\input{Content/Abstract}
\end{abstract}

\maketitle

\section{Introduction}\label{sec:Intro}
\input{Content/Intro}

\section{Preliminaries}\label{sec:Preli}
\input{Content/Preli}

\section{Rotational LW surfaces}\label{sec:rotLW}
\input{Content/RotLW}

\section{Rotational cmc surfaces}\label{sec:channelCMC}
\input{Content/RotCMC}

\section{Classification of channel linear Weingarten surfaces in 
  \texorpdfstring{$\H^3$}{H3}}\label{sec:channelLW}
\input{Content/ChannelLW}

\section*{Acknowledgements}
The author would like to express his gratitude to Udo Hertrich-Jeromin for many
years of mentoring and suggesting interesting avenues of research, like this 
one. Further, the author thanks Joseph Cho, Yuta Ogata, Mason Pember, and Wayne 
Rossman for many fruitful discussions about this subject. This work 
was done while the author was a JSPS International Research Fellow (Graduate 
School of Science, Kobe University) and has been supported by the JSPS 
Grant-in-Aid for JSPS Fellows 22F22701.

\printbibliography

\begin{minipage}[t][3cm][c]{0.6\columnwidth}
  \textbf{Denis Polly}\\
  School of Computation, Information and Technology\\
  TU M\"{u}nchen\\
  Boltzmannstra{\ss}e 3/3\\
  85748 Garching bei M\"{u}nchen\\
  \url{denis.polly@tum.de}
\end{minipage}

\end{document}

%% file: Preamble/Environments.tex
  \theoremstyle{remark} 
  \newtheorem{defi}{Definition}[section]
	\crefname{defi}{definition}{definitions}
	
  \newtheorem{bsp}[defi]{Example}
  
  \newtheorem{examples}[defi]{Examples}

  \newtheorem{rem}[defi]{Remark}
	\crefname{rem}{remark}{remarks}
  
  \newtheorem{rems}[defi]{Remarks}
  \crefname{rems}{Remarks}{remarks}

\theoremstyle{plain} 
  \newtheorem{thm}[defi]{Theorem}
	\crefname{thm}{theorem}{theorems}
  
  	\newtheorem{prop}[defi]{Proposition}
	\crefname{prop}{proposition}{propositions}

    \newtheorem{lem}[defi]{Lemma}
	\crefname{lem}{lemma}{lemmata}

	\crefname{cor}{corollary}{corollaries}

%% file: Preamble/Commands.tex
\newcommand{\N}{\mathbb{N}}
\newcommand{\Q}{\mathbb{Q}}

\newcommand{\R}{\mathbb{R}}

\newcommand{\newcal}[2]{\newcommand{#1}{\mathcal{#2}}}
\newcommand{\renewcal}[2]{\renewcommand{#1}{\mathcal{#2}}}
\newcommand{\newfrak}[2]{\newcommand{#1}{\mathfrak{#2}}}
\newcommand{\renewfrak}[2]{\renewcommand{#1}{\mathfrak{#2}}}

\newcommand{\renewbb}[2]{\renewcommand{#1}{\mathbb{#2}}}
\renewbb{\H}{H}
\renewbb{\S}{S}
\newfrak{\s}{s}
\newfrak{\p}{p}
\newfrak{\q}{q}
\newfrak{\e}{e}
\renewfrak{\v}{v}
\newfrak{\w}{w}
\renewfrak{\o}{o}
\newfrak{\f}{f}
\newfrak{\n}{n}
\renewfrak{\c}{c}
\renewcal{\L}{L}
\newcal{\Lines}{Z}
\renewcommand{\Q}[1]{\mathfrak{R}_{#1}}
\newcommand{\T}[1]{\mathfrak{T}_{#1}}
\newcommand{\Iso}[1]{\operatorname{Iso}_{#1}}

\newcommand{\jac}[2]{\operatorname{#1}_{#2}} 
\newcommand{\Pii}[3]{\Pi_{#1}\!\left(
  #2; #3
\right)}
\newcommand{\ro}[2]{\rho^{#1}_{#2}}
\newcommand{\tra}{\top}
\newcommand{\sgn}{\operatorname{sgn}}
\newcommand{\FAC}{\Xi}
\newcommand{\spn}[1]{\left\langle#1\right\rangle}

	\definecolor{newred}{RGB}{113, 21, 33}
	\definecolor{newblue}{RGB}{93, 105, 112}
  \definecolor{newblack}{RGB}{10,10,10}
  \definecolor{newwhite}{RGB}{240,240,240}

%% file: Content/Abstract.tex
We give a classification of rotational cmc surfaces in non-Euclidean space forms
in terms of explicit parametrizations using Jacobi elliptic functions. Our 
method hinges on a Lie sphere geometric description of rotational linear 
Weingarten surfaces and can thus be applied to a more general class of surfaces.
As another application of this framework, we give explicit parametrizations of a
class of rotational constant harmonic mean curvature surfaces in hyperbolic 
space. In doing so, we close the last gaps in the classification of all
channel linear Weingarten surfaces in space forms, started in 
\cite{hertrich-jeromin2023}.

%% file: Content/Intro.tex
The study of surfaces of constant mean curvature (cmc surfaces for short) has 
been a subject of interest in surface theory since the 19\textsuperscript{th} 
century. It was Delaunay \cite{delaunay1841} who
discovered that the profile curves of cmc surfaces of revolution, today 
called \emph{Delaunay surfaces}, can be constructed by rolling conic sections 
along the axis of revolution and tracing their focal points. This alludes to a 
connection between the profile curves of Delaunay surfaces
and conics, also on the level that Jacobi elliptic differential equations
can be used to describe them. This connection 
extends to the wider class of \emph{linear Weingarten surfaces (lW)} of 
revolution, that is, surfaces with an affine relationship between their Gauss
and mean curvature \cite{hertrich-jeromin2015}.

Surfaces invariant under rotations with constant mean curvature
have been of interest in many contexts, extending to space forms of constant 
curvature $\kappa \neq 0$. While rotational surfaces in 
spherical space forms ($\kappa>0$) are invariant under a $1$-parameter family of
Euclidean rotations in $\R^4$, the situation is much richer in cases in 
hyperbolic space forms ($\kappa<0$), where rotations come in three different 
flavors \cite{carmo1983}. 

The situation is even richer in hyperbolic space forms, given that the family of
surfaces of constant mean curvature $H$ splits into three distinct classes 
depending on the sign of $H^2 + \kappa$: surfaces with $H^2 + \kappa \geq 0$ 
have been subject to many publications, being closely 
related to cmc (minimal) surfaces in $\R^3$ and thus allowing for investigation 
via representations like the DPW method \cite{dorfmeister2014} or the Bryant 
representation \cite{bryant1987}. There is not as much literature in the 
$H^2+\kappa < 0$ case, however some publications investigate general cmc 
surfaces, hence apply also to this case 
\cite{abe2018,castillon1998,dorfmeister2014} and \cite{mori1983}.

The study of rotational cmc surfaces has come back into focus recently, providing
multiple classification theorems in hyperbolic and spherical space forms from
differing points of view. This includes description of the profile curves of
these surfaces as energy minimizers \cite{arroyo2019}, as level sets 
of an auxiliary potential \cite{barros2012} or their explicit parametrization in
 terms of one of their principal curvatures \cite{dursun2020}. 

\begin{figure}
\centering
  \begin{tikzpicture}
    \node at (-5,6) {\includegraphics[width=0.3\textwidth]{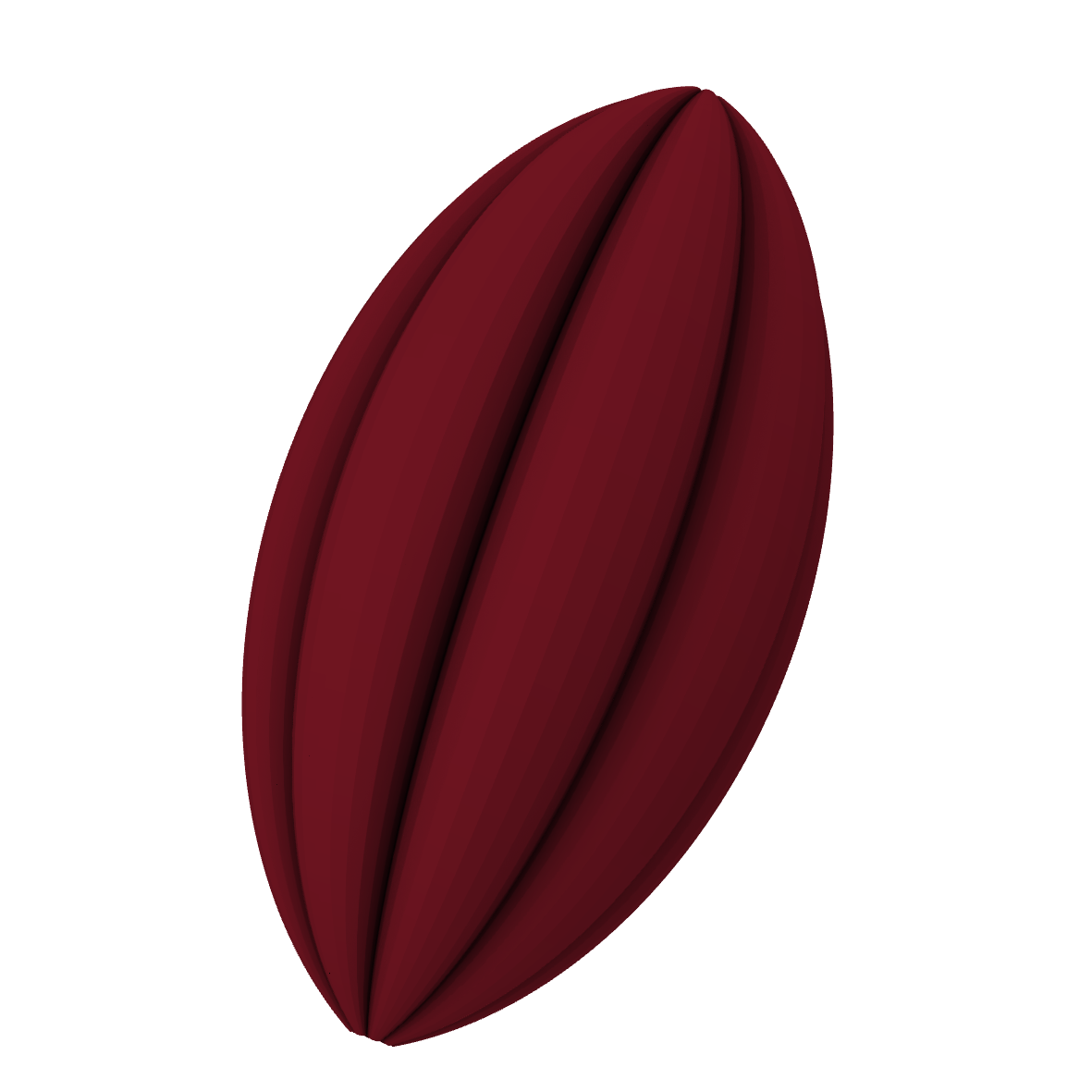}};
    \node at (0,6) {\includegraphics[width=0.3\textwidth]{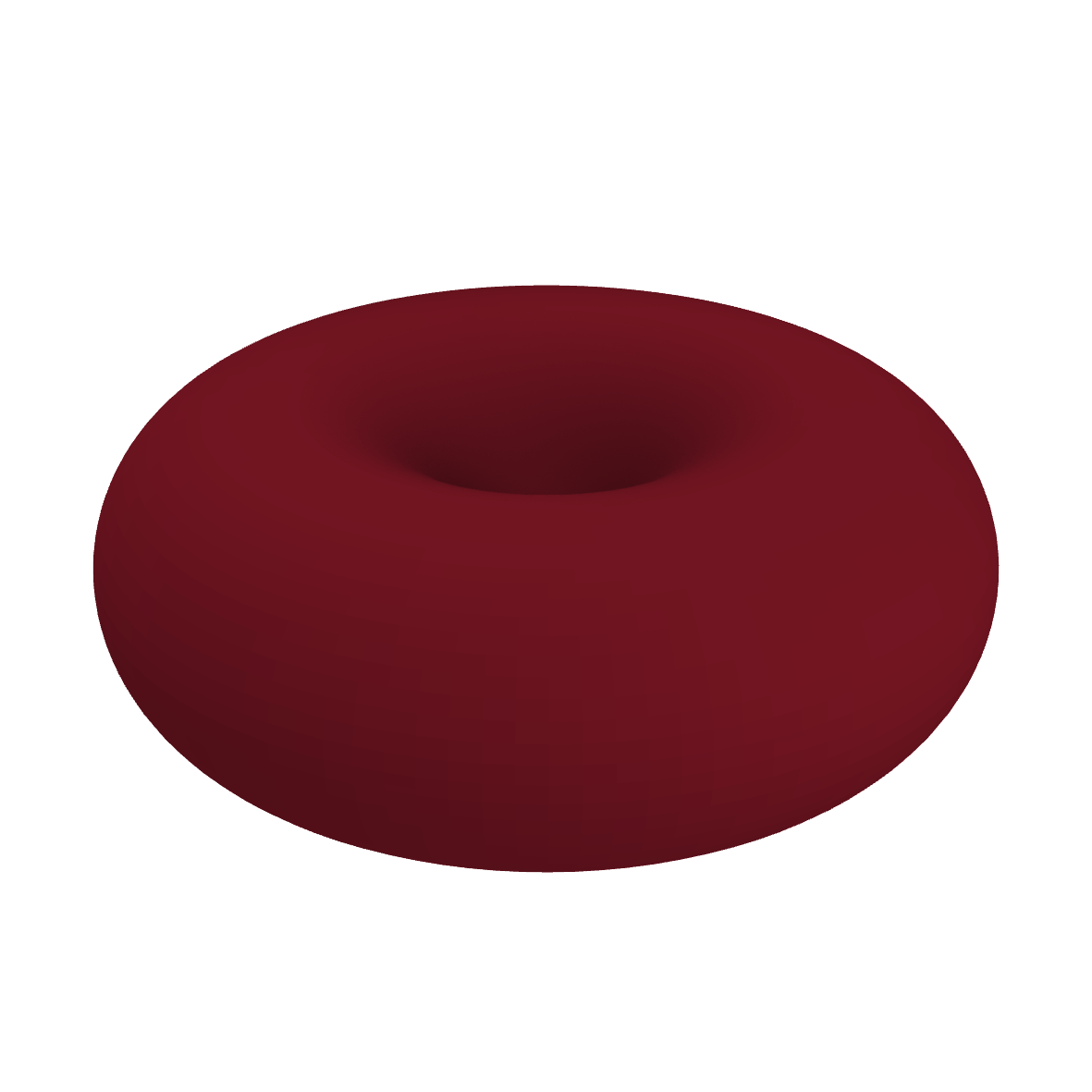}};
    \node at (5,6) {\includegraphics[width=0.3\textwidth]{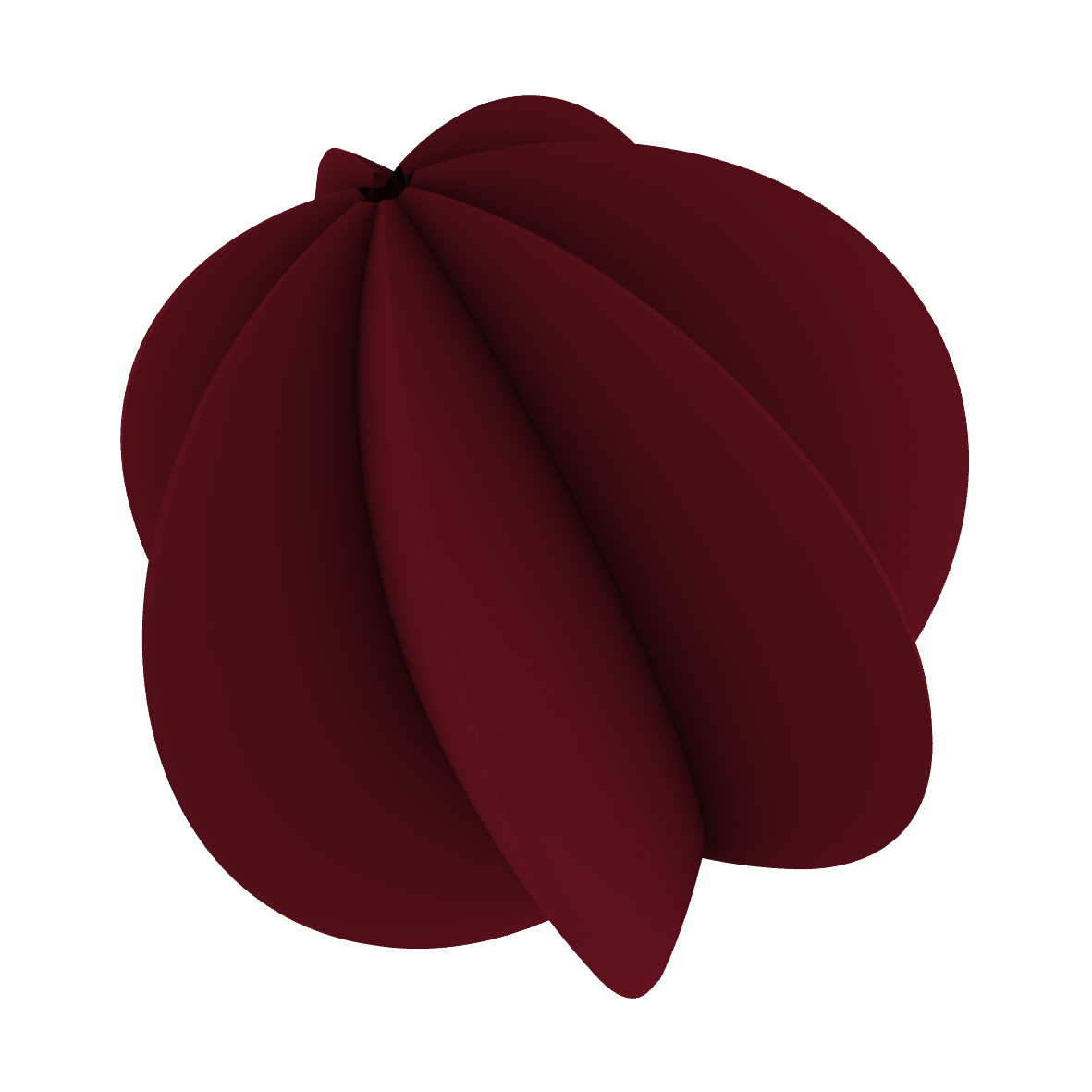}};
    \node at (-3,1) {\includegraphics[width=0.3\textwidth]{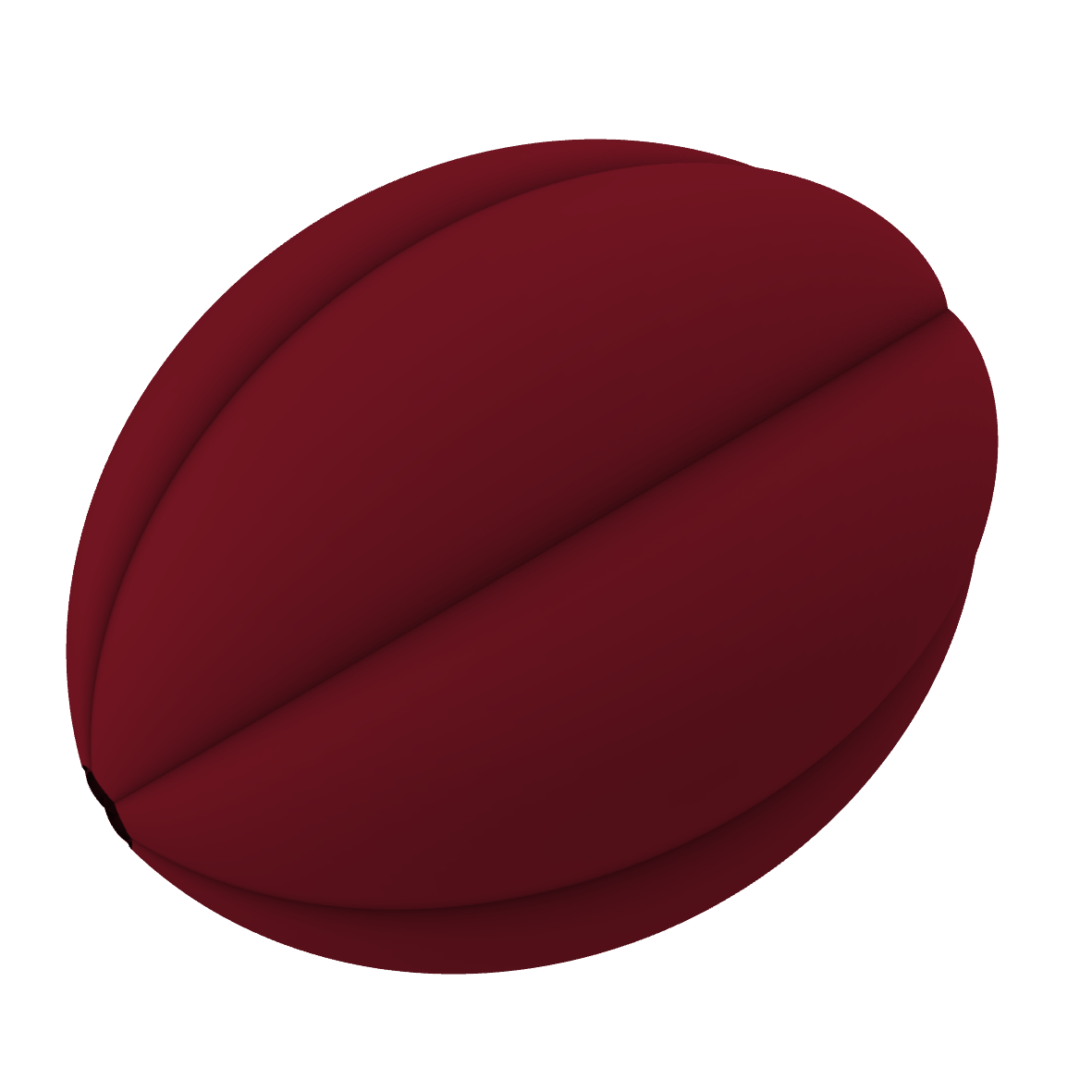}};
    \node at (3,1) {\includegraphics[width=0.3\textwidth]{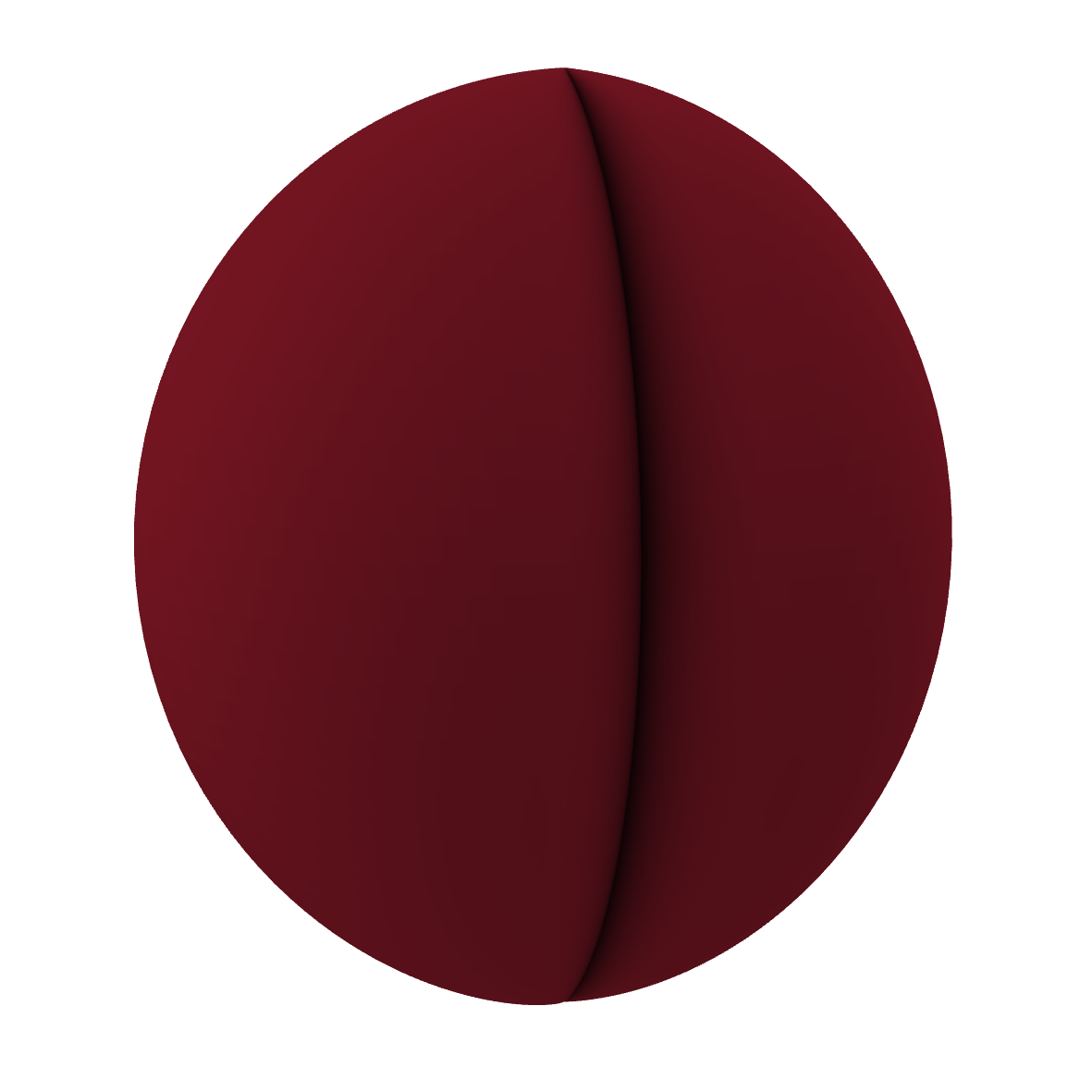}};
    \node at (-5,3.4) {a)};
    \node at (0,3.4) {b)}; 
    \node at (5,3.4) {c)};
    \node at (-3,-1.5) {d)};
    \node at (3,-1.5) {e)};
  \end{tikzpicture}
  \caption{Fruits in the linear Weingarten: a) cocoa pot surface (constant Gauss
    curvature in $\H^3$, hyperbolic rotation) b) saturn peach surface (cmc 
    in $\S^3$) c) star fruit surface (constant Gauss curvature in $\H^3$, 
    hyperbolic rotation) d) cantaloupe surface (cmc in $\H^3$, hyperbolic 
    rotation) e) peach front (intrinsically flat in $\H^3$, hyperbolic 
    rotation). Plotted and rendered using Rhino and mathematica.
    }
  \label{fig:Fruits}
\end{figure}
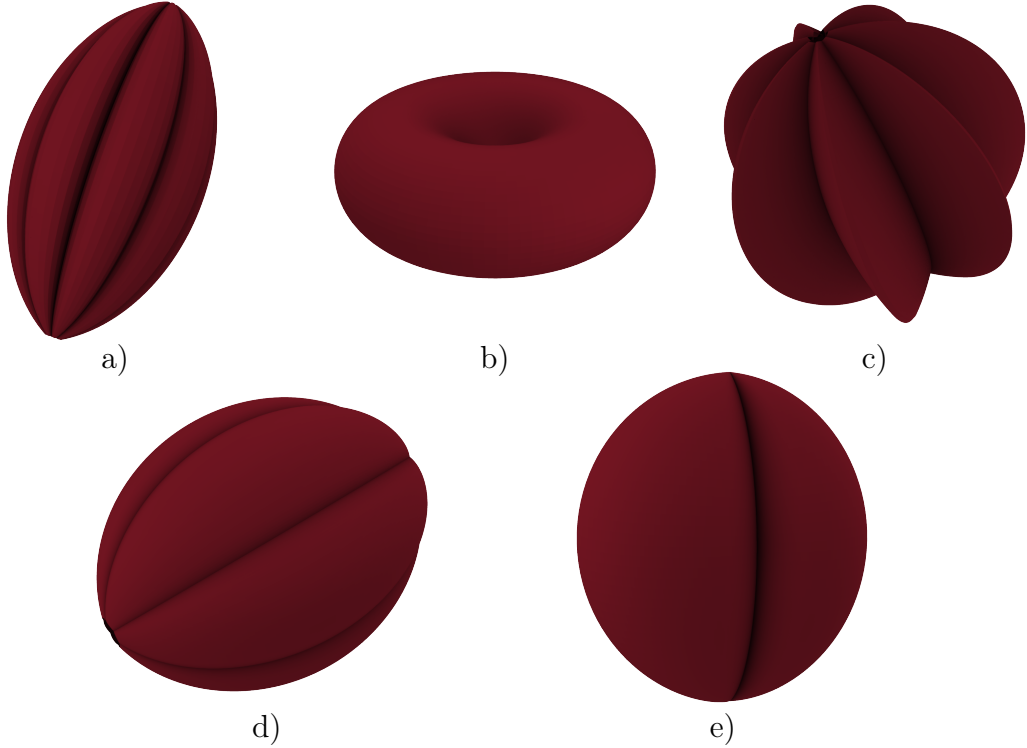

Some of the references just mentioned consider the wider class of rotational
lW surfaces in space forms. This class was 
considered in \cite{hertrich-jeromin2023}, where the authors proved that
every lW surface enveloping a $1$-parameter family of spheres is rotational 
in its space form. Furthermore, the authors shed light on the connection
between rotational lW surfaces and Jacobi elliptic functions by presenting
explicit parametrizations of a wide class of rotational lW surfaces in terms
of these functions, namely those that are parallel to surfaces of constant 
Gauss curvautre (cGc), see \Cref{fig:Fruits} a,c and e. However, the specific 
case of surfaces parallel to cmc $H$ surfaces in hyperbolic space forms with 
$H^2 + \kappa<0$ was unresolved: these surfaces are not related to a constant 
Gauss curvature surface via parallel transformations, hence, no parametrizations
for them can be obtained from results stated in \cite{hertrich-jeromin2023}.

The main motivator of this paper is to close this gap in the classification
of channel lW surfaces in space forms via Jacobi elliptic functions. We will 
develop a new approach to derive differential equations governing
the profile curves of rotational lW surfaces in non-Euclidean space forms. This
new approach stems from a characterization of lW surfaces via enveloped 
isothermic sphere congruences with additional properties \cite{burstall2012}. 
This description belongs to the realm of Lie sphere geometry and will prove
particularly useful once applied to the specific case of cmc surfaces. Further,
the following advantages of our setup shall be emphasized: 
\begin{itemize}
  \item since spherical and hyperbolic space form geometries appear as 
  subgeometries of Lie sphere geometry, we can carry out our analysis for
  both ambient spaces simultaneously, breaking symmetry only at the 
  very end to obtain space form specific parametrizations. 
  \item surfaces of constant harmonic mean curvature can be investigated
  with similar ease. This is useful because another
  class of surfaces, not covered in the classification results of 
  \cite{hertrich-jeromin2023}, are surfaces with constant harmonic mean 
  curvature greater than or equal to $1$. 
\end{itemize}
The basic idea to derive parametrizations of rotational cmc surfaces in space 
forms was laid out in the authors PhD thesis \cite{polly2022}. Our goal is to 
bring this idea to fruition. 

This paper is organized as follows: in Section~\ref{sec:Preli} we set the stage 
of Lie sphere geometry and its space form subgeometries. We will discuss briefly
how lW and rotational surfaces can be investigated using this setup but refer 
the details-seeking reader to publications that discuss this setup in more 
detail, such as \cite{burstall2014, cecil2008,hertrich-jeromin2023}. This 
section will also introduce a characterization of lW surfaces in terms of 
isothermic sphere congruences that take values in certain families of spheres
\cite{burstall2012} and are, in some sense, 
Christoffel dual (\Cref{lem:ChristoffelDualLifts}). In Section~\ref{sec:rotLW} 
we will further investigate these isothermic sphere congruences in the case
of rotational lW surfaces in non-Euclidean space forms, leading to differential 
equations governing their profile curves in \Cref{prop:JEF,prop:parJEF}. 
Section~\ref{sec:channelCMC} will see us applying these differential equations 
to the case of rotational cmc surfaces. This is particularly convenient as it
will turn out that one of the isothermic sphere congruences in this case is 
the surface itself. Thus, we arrive at explicit parametrisations of rotational
cmc surfaces in non-Euclidean spaceforms. These results split along the lines
of ambient space form, type of rotation and the sign of $H^2+\kappa$
(\Cref{thm:DelaunayS3,thm:DelaunayH3,thm:parDelaunayH3,thm:smallH,%
thm:smallHpar}). 

Finally, in Section~\ref{sec:channelLW}, we finish the 
classification of channel lW surfaces in terms of Jacobi elliptic functions, 
started in \cite{hertrich-jeromin2023}. We show how the formulas 
derived in Seciton~\ref{sec:rotLW} can be applied to the case of surfaces with 
constant harmonic mean curvature, which will allow us to close the gap in the 
classification of rotational lW surfaces in hyperbolic space forms 
(\Cref{thm:chc,thm:cLWH3}).

%% file: Content/Preli.tex
In this section, we will introduce the light cone model of Lie sphere geometry
and explain how to investigate linear Weingarten surfaces, channel surfaces and
rotationally invariant surfaces in this set-up. Our description of linear 
Weingarten surfaces in Lie sphere geometry was originally given in 
\cite{burstall2012}. 

We are interested in rotational cmc surfaces. However, we will discuss the 
more general class of channel linear Weingarten surfaces in preparation
of the special cases dealt with in the following sections. 
A more detailed account of these preparations can 
be found in \cite{hertrich-jeromin2023}, where a Lie sphere geometric 
investigation of channel linear Weingarten surfaces was laid out. In this note
we will state the necessary results and refer to \cite{burstall2012} and 
\cite{hertrich-jeromin2023} for proofs. 

\subsection{Lie sphere geometry}\label{sec:LSG}
\input{Content/Preli/LSG}

\subsection{Channel linear Weingarten surfaces}\label{sec:CLW}
\input{Content/Preli/CLW}

\subsection{Rotational surfaces}\label{sec:rotationalSurfaces}
\input{Content/Preli/Rot}

%% file: Content/Preli/LSG.tex
Let $\R^{4,2}$ be the $6$-dimensional real vector space with non-degenerate 
inner product $(\cdot, \cdot)$ of signature $(++++--)$. 
A vector $\v\in\R^{p,q}$ is called 
\emph{timelike, spacelike} or \emph{lightlike} if $(\v,\v)$ is negative, 
positive or vanishes respectively. We denote by $L^{4,1}$ the set
of all lightlike vectors (the \emph{light cone}).

We call a basis of $\R^{4,2}$ 
\emph{orthonormal} if the inner product takes the form
\begin{align*}
  (\v, \w) = -v_0w_0 + v_1 w_1 + v_2 w_2
    + v_3 w_3 + v_4 w_4 - v_5 w_5,
\end{align*}
and \emph{pseudo-orthonormal}, if the inner product takes the form
\begin{align*}
  (\v, \w) = -v_0w_5 - v_5w_0 + v_1 w_1 
    + v_2 w_2 + v_3 w_3 + v_4 w_4.
\end{align*}

We define the set of \emph{(L-)spheres} in a space form as the union of the sets
 of points, oriented hyperspheres and oriented hyperplanes (i.e. complete, 
totally geodesic hypersurfaces). In the light cone model of Lie sphere geometry,
originally introduced by Lie in \cite{lie1872} (see \cite{cecil2008} for a 
modern account), L-spheres of $3$-dimensional space forms are represented by 
points in the projective light cone 
\begin{align*}
  \L=\mathbb{P}L^{4,1} = \{\spn{\v} \subset \R^{4,2}: \v\neq 0 
  \textrm{ is lightlike}\}.
\end{align*}
which is called the \emph{Lie quadric}. Two L-spheres $s_1, s_2 \in \L$ are in 
\emph{oriented contact} if and only if there are homogeneous coordaintes 
$\s_i \in s_i$--which we also call \emph{lifts of $s_i$}--such that 
\begin{align*}
  (\s_1, \s_2) = 0.
\end{align*}
Note that this equation is independent of the choice of lifts. 
We call a maps $s:\Sigma \to \L$ into the Lie quadric \emph{sphere congruences}.

The transformation group of Lie sphere geometry
is derived from the the projective orthogonal group of $\R^{4,2}$: 
$O(4,2)/\{\pm \operatorname{Id}\}$ acts on the Lie quadric in the way linear 
transformations typically act on projective spaces (see 
\cite[Cor 3.3]{cecil2008}). It maps oriented spheres to oriented spheres and 
preserves oriented contact. Further, it preserves \emph{linear sphere 
complexes}: each vector $\p \in \R^{4,2}$ defines a family of spheres  
\begin{align*}
  \{s\in\L: (\p, \s) = 0~\forall \s \in s\}
\end{align*}
called the \emph{(linear) sphere complex defined by $\p$}. We also
call $\p$ itself a (linear) sphere complex. Note that two parallel vectors
define the same sphere complex. We call a sphere complex 
\emph{timelike, lightlike} and \emph{spacelike} if the spanning vector $\p$
has the respective causal character. 

We may break symmetry and recover models for the space form geometries: let 
$\p$ be a unit timelike sphere complex and $\q$ a perpendicular sphere
complex. The affine sub-quadric
\begin{align*}
  \Q{\p, \q} = \{\v \in L^{4,1}: (\v,\p) = 0, (\v,\q) =-1 
  \},
\end{align*}
has constant sectional curvature $\kappa = -(\q, \q)$ 
(see \cite[Sec 2]{burstall2018}). Since we define 
L-spheres in the sphere complex $\p$ to be points in the space form, 
we call $\p$ a \emph{point sphere complex}. Similarly,  
$\q$ defines the curvature of the space form, hence, we call it the \emph{space 
form vector}. For a point sphere $s$ we call $\s \in \Q{\p,\q}$ with 
$\spn{\s}=s$ a \emph{space form lift}. Together with the subgroup 
$\Iso{\p, \q}$ of Lie sphere transformations fixing $\p$ and $\q$, $\Q{\p,\q}$ 
is a model for a space form geometry, see \cite[Sect 2]{burstall2018} or 
\cite[Chap 2]{cecil2008}. The set of oriented hyperplanes in $\Q{\p, \q}$ is
 represented by points in the affine sub-quadric 
\begin{align*}
  \T{\p, \q} = \{\v \in L^{4,1}: (\v,\p)=-1, (\v,\q)=0
  \}.
\end{align*} 
Sphere congruences taking values in the sphere complex $\q$ are thus called 
\emph{plane congruences}.

The study of hypersurfaces in space forms is carried out in Lie sphere geometry
via the utilization of \emph{Legendre lifts}: let $\f:\Sigma \to \Q{\p, \q}$ 
be a parametrization of a  hypersurface in the space form $\Q{\p,\q}$ and 
$\n:\Sigma \to \T{\p, \q}$ 
its \emph{tangent plane congruence}. Then $\f$ and $\n$ satisfy
\begin{align}\label{eq:contactConditions}
  (\f, \n) = 0, \textrm{ and } (d\f, \n) = 0.
\end{align}
This implies that the (projective) line congruence $\Lambda = \spn{\f, \n}$ is a
\emph{Legendre immersion} (see \cite[Sect 4.2]{cecil2008}).  We call 
\eqref{eq:contactConditions} the \emph{contact conditions}. Additionally we 
assume the \emph{immersion condition}
\begin{align*}
  \forall p\in\Sigma\, \forall \s \in \Gamma\Lambda: d_p\s(X)\in\Lambda(p)
  \Rightarrow X=0.
\end{align*}
We say that $\Lambda$ \emph{envelopes a sphere congruence $s=\spn{\s}$} 
if $\s(p) \in \Lambda(p)$ for all $p \in \Sigma$. We will denote the set of 
sphere congruences enveloped by $\Lambda$ by $\Gamma \Lambda$ (similarily, we 
shall denote the set of maps into $L^{4,1}$ such that $\s(p) \in s(p)$ for 
a given sphere congruence as $\Gamma s$)\footnote{
This notation is motivated by the interpretation of $\Lambda$ (or $s$) as 
a sub-bundle of the trivial vector bundle $\Sigma^2\times \R^{4,2}$. Then,
the set of rank $1$ subbundles of $\Lambda$ (sections of $s$)--typically denoted
$\Gamma\Lambda$ ($\Gamma s$)--is the set of enveloped sphere congruences (lifts).
}.

Lie sphere transformations naturally act on Legendre lifts. Together with 
the following lemma, this implies that Lie sphere transformations act on 
surfaces in space forms in a meaningful way. 

\begin{lem}[{\cite[Sec 2.5]{cecil2008}}]
Given a point sphere complex $\p$ a Legendre immersion $\Lambda$, there is 
precisely one point sphere congruence $f$ enveloped by $\Lambda$. 
\end{lem}

Similarly, any Legendre immersion $\Lambda$ generically envelopes precisely one 
(tangent) plane congruence\footnote{We exclude 
the degenerate case where $\Lambda$ consists of parallel planes in a hyperbolic 
space form.} for any space form vector $\q$. Additionally, 
$\Lambda$ envelopes \emph{curvature sphere congruences} 
\cite[Sec 4.4]{cecil2008}: there are coordinates $(u,v)$ and sphere congruences 
$s_1, s_2 \in \Gamma \Lambda$ such that 
\begin{align*}
  (\s_1)_u, (\s_2)_v \in \Gamma\Lambda,
\end{align*}
for any lifts $\s_i \in s_i, ~i=1,2$. The coordinates $(u,v)$ are then curvature
line coordinates for any space form lift of any point sphere envelope. Away from
umbilics we have $\Lambda = s_1 \oplus s_2$.

%% file: Content/Preli/CLW.tex
This subsection is a brief summary of the Lie sphere geometric treatment of 
channel linear Weingarten surfaces given in \cite{hertrich-jeromin2023}. For the
purposes of this paper, we present \Cref{thm:ClWInLSG} and \Cref{thm:CLWIsRot} 
and point to \cite{hertrich-jeromin2023} and references therein for details. We 
conclude with a short prove of a lemma that will be instrumental in the following
investigations.

Given a point sphere complex $\p$ and a space form vector $\q$, 
the curvature of an L-sphere in $\Q{\p, \q}$ can be computed via
\begin{align*}
  \kappa_s = \frac{(\s, \q)}{(\s, \p)},
\end{align*}
which is lift-independent \cite{burstall2012}. We can apply this to recover the 
principal curvatures $\kappa_1, \kappa_2$ of a surface $\f$ with values in 
$\Q{\p, \q}$ from the curvature sphere congruences $s_1, s_2$ of its Legendre 
lift $\Lambda$. We define the \emph{mean curvature $H$} and the 
\emph{(extrinsic) Gauss curvature} $K$ of $\f$ via
\begin{align*}
  H=\frac{\kappa_1 + \kappa_2}{2}, \quad K = \kappa_1\kappa_2.
\end{align*}
A surface parametrized by $\f$ is \emph{linear Weingarten (lW)} if there exists 
a (non-trivial) triple of constants $a, b, c \in \R$ such that
\begin{align}\label{eq:LW}
  aK + 2bH + c = 0.
\end{align}
The Discriminant $D$ of an lW surface is defined as 
$D = b^2 - ac$ and we call the surface \emph{tubular} if $D =0$. The class of
lW surfaces contains important sub-classes: for $a=0$ we obtain \emph{constant 
mean curvature (cmc) surfaces}, for $b=0$ \emph{constant Gauss curvature (cGc) 
surfaces}. For $c=0$ we learn that the \emph{harmonic mean curvature 
$\overline{H} := K/H$} is constant; we call such surfaces \emph{constant 
harmonic mean curvature (chc) surfaces}. 

In \cite{burstall2012}, the authors prove that linear Weingarten surfaces
in space forms can be characterized in Lie sphere geometry as a subclass of 
the class of $\Omega$-surfaces. We will state this result and point to 
the original source for a proof. 

We call a sphere congruence $s$ \emph{isothermic} if it
possesses a Moutard lift, that is,
\begin{align*}
  \exists \s \in \Gamma s: \s_{uv} || \s,
\end{align*}
for suitable coordinates $(u,v)$ which are then curvature line
parameters. Following Demoulin 
\cite{demoulin1911, demoulin1911a}, we call a Legendre immersion an
\emph{$\Omega$-surface} if it envelopes a (possibly complex conjugated)
pair of isothermic sphere congruences that separate the curvature sphere 
congruences harmonically. We also call a surface $\f$ in a space form $\Omega$
if its Legendre lift is an $\Omega$-surface.

\begin{thm}[{\cite{burstall2012}}]\label{thm:ClWInLSG}
  The Legendre lifts of non-tubular linear Weingarten surfaces are those
  $\Omega$-surfaces $\Lambda= \spn{\s^+, \s^-}$ whose isothermic sphere 
  congruences $s^\pm$ take values in fixed linear sphere complexes $\p^\pm$. 
  The plane $\spn{\p^+, \p^-}$ is spanned by the point sphere complex $\p$
  and the space form vector $\q$.
\end{thm} 

\begin{rem}
  The isothermic sphere congruences $s^\pm$ of a linear Weingarten surface 
  are complex conjugate if and only if $\Delta<0$. 
\end{rem}

We therefore generally call an $\Omega$-surface $\Lambda$ a \emph{linear
Weingarten surface} if its isohtermic sphere congruences take values in 
fixed linear sphere complexes $\p^\pm$. To emphasize this, we denote the lW
surface by $(\Lambda, \p^\pm)$. 

\begin{bsp}\label{bsp:conservedQuantities}
  Let $\f:\Sigma \to \Q{\p, \q}$ be a parametrization of a cmc $H$ surface. 
  According to \cite[Sec 4.2]{burstall2018}, the sphere complexes $\p^\pm$ are 
  given by
  \begin{align*}
    \p^+ = \p, \quad \p^- = \q - H\p. 
  \end{align*}
  Since the point sphere envelope of a Legendre immersion is unique, this implies
  that the surface itself is one of the isothermic sphere congruences spanning 
  $(\Lambda, \p^\pm)$, $s^+ = \spn{\f}$. The other one is given by
  \begin{align*}
    \s^- = \n + H\f.
  \end{align*}
  This is the sphere congruence enveloped by $\Lambda$ that has in each point 
  the same mean curvature as the surface. According to 
  \cite[\textsection 67]{blaschke1929} this is the central sphere congruence of
  the surface (or its conformal Gauss map, as introduced in
  \cite{bryant1984}). Further, $\s^\pm$, with $\s^+ = \f$ 
  are Christoffel dual lifts of $s^\pm$ (see \Cref{lem:ChristoffelDualLifts} 
  below). 
  
  Somewhat dual to this is the description of chc surfaces: let 
  $\f:\Sigma \to \Q{\p,\q}$ be a parametrization of a chc $\overline{H}$ 
  surface. Then we have
  \begin{align*}
    &\p^+ = \q, \quad \p^- = \p - \frac{1}{\overline{H}}\q \\
    \Rightarrow ~ &\s^+ = \n, \quad \s^- =\f + \frac{1}{\overline{H}}\n.
  \end{align*}
  This means that one isothermic sphere congruence is given by the tangent
  plane congruence and the other is the sphere congruence that has the same 
  harmonic mean curvature as the surface (this is called the middle sphere
  congruence in \cite[\textsection 67]{blaschke1929}
  or Laguerre Gauss map in \cite{musso2019}). 
\end{bsp}

Next, consider a channel lW surface, that is, the surface envelopes a 
$1$-parameter family of spheres. The enveloped sphere curve $s$ is then
automatically a curvature sphere congruence (as it has a lift that only
depends on one curvature line parameter). The simplest example is given by
rotational surfaces which we shall describe in the next subsection. 

The following theorem has been shown in \cite{hertrich-jeromin2023}, where we 
point to for a proof. 

\begin{thm}[{\cite[Thm 4.7]{hertrich-jeromin2023}}]\label{thm:CLWIsRot}
  Every non-tubular channel linear Weingarten surface in a space form 
  $\Q{\p,\q}$ is a rotational surface. 
\end{thm}

By this theorem, all channel cmc surfaces are rotational in their space form.  
Further, because $D = b^2 > 0$ for cmc surfaces, they envelope a real pair
of isothermic sphere congruences. Thus, it will be our goal to describe  
rotational $\Omega$ surfaces $(\Lambda, \p^\pm)$ with $\p^\pm$ a real pair. 

We close this section with the following lemma, introducing a property 
of the isothermic sphere congruences enveloped by an $\Omega$-surface 
that is important for our later investigations. 

\begin{lem}[{\cite[Subsec 4.2.2]{pember2015}}]\label{lem:ChristoffelDualLifts}
  Let $\Lambda$ be an $\Omega$-surface spanned by two isothermic sphere 
  congruences $s^\pm$ and let $(u,v)$ be curvature line coordinates. Then there 
  are lifts $\s^+\in\Gamma s^+$ and $\s^-\in \Gamma s^-$ such that
  \begin{align}\label{eq:ChristoffelDual}
    \begin{split}
      \s^-_u &= CU\frac{\s^+_u}{|\s^+_u|^2} \\
      \s^-_v &= -CV\frac{\s^+_v}{|\s^+_v|^2},
    \end{split}
  \end{align}
  where $U, V$ are functions of (only) $u$ and $v$ respectively. These lifts are
  unique up to constant rescaling. We call them the \emph{Christoffel dual 
  lifts}. 
\end{lem}

\begin{proof}
  We will only sketch the proof of this lemma, for details see 
  \cite[Lem 3.6]{polly2022}\footnote{
  The cited result only concerns the case of real isothermic sphere congruences.
  However, with the specific lifts given in this proof, the computations can be
  generalized to include the complex conjugate case.}
  The existence of Christoffel dual lifts follows from the following 
  characterization of $\Omega$ surfaces: there exist lifts 
  $\s_i \in s_i$ of their curvature sphere congruences such that
  \begin{align*}
    (\s_1)_u=\phi_u \s_2, \quad (\s_2)_v = \varepsilon^2 \phi_v \s_1,
  \end{align*}
  where $\varepsilon \in \{i, 1\}$, see \cite[Prop 2.4]{polly2022}. 
  With these special lifts it is straightforward to prove that the lifts 
  \begin{align*}
    \c^\pm = \overline{e^{\mp \varepsilon \phi}}(\s_1 \pm \varepsilon \s_2),
  \end{align*}
  are Christoffel dual lifts of $s^\pm$. 
  
  Now, assume there is another pair of Christoffel dual lifts, given by 
  $\alpha^\pm \c^\pm$. Checking \eqref{eq:ChristoffelDual} then proves that 
  $\alpha^\pm$ have to be constants.
\end{proof}

\begin{rem}\label{rem:ChristoffelDual}
  The functions $U, V$ have the property that
  \begin{align*}
    \frac{(\c^+_u, \c^+_u)}{U} = \frac{(\c^+_v, \c^+_v)}{V}.
  \end{align*}
  The existence of functions with these properties follows from the existence
  of a Moutard lift \cite[Section 2.3]{polly2022}.
\end{rem}

%% file: Content/Preli/Rot.tex
The goal of this subsection is to give a quick overview over the 
unifying description of rotational surfaces in space forms within the
realm of Lie sphere geometry. We will shoot for brevity and refer the 
interested reader to \cite[Sec 3]{hertrich-jeromin2023}. 

As we discussed, the subgroup $\Iso{\p, \q}$ of Lie sphere transformations 
that fix a point sphere complex $\p$ and a space form vector $\q$ models
the isometry group of the space form $\Q{\p, \q}$.
Let $\Pi$ be a $2$-plane with $\Pi \perp \spn{\p, \q}$, we 
call a $1$-parameter subgroup $\rho$ of $\Iso{\p, \q}$ that acts as the 
identity on $\Pi^\perp$, a \emph{subgroup of rotations}. We call $\rho$
\begin{itemize}
	\item \emph{ellitpic} if $\Pi$ has signature $(++)$
  \item \emph{parabolic} if $\Pi$ has signature $(+0)$
  \item \emph{hyperbolic} if $\Pi$ has signature $(+-)$.
\end{itemize}
Note that, since $\p$ is timelike, these are the only signatures that can 
occur and that not all types of rotations exist in $\S^3$ and $\R^3$. We will
utilize a basis $\{\v, \e\}$ of $\Pi$ that is orthogonal with $\e$ being unit 
spacelike and $(\v, \v)$ encoding the signature.

\begin{defi}
A sphere congruence $s:\Sigma^2 \to \L$ is called \emph{rotational} if it can be
parametrized as\footnote{The subgroup of rotations acts on the projective curve
$c$ via $\ro{}{}c :=\spn{\ro{}{}\c}$ for any lift $\c \in \Gamma c$.} 
\begin{align*}
  s(\theta, t) := \ro{}{\theta} c(t),
\end{align*}
where $\ro{}{}$ is a $1$-parameter subgroup of rotations in a space form 
$\Q{\p, \q}$ with respect to $\Pi \perp \spn{\p,\q}$ and 
$c:=s(0,\cdot):I \to \spn{\e}^\perp\cap \L$. We call $c$ its \emph{planar 
profile curve}. A Legendre immersion $\Lambda$ is called \emph{rotational} with
respect to $\ro{}{}$ if it is spanned by two rotational sphere 
congruences.
\end{defi}

As a special case, a surface $\f:\Sigma^2 \to \Q{\p, \q}$ is called rotational 
if its lift $f=\spn{\f}$ is rotational. Note that not any sphere congruence 
enveloped by a rotational Legendre immersion is rotational. However, in the case
of Legendre lifts, some prominent envelopes are. 

\begin{lem}\label{lem:rotationalLegendreLift}
  The Legendre lift $\Lambda$ of a rotational surface $f$ and its curvature 
  sphere congruences are rotational. If $f$ is an $\Omega$-surface, then the
  enveloped isothermic sphere congruences $s^\pm$ are rotational as well.
\end{lem}

\begin{proof}
  Let $\c \in \Gamma c$ be the space form lift of the profile curve of $f$. 
  There is a unique curve ${\n_0: \Sigma^2 \to \spn{\e}^\perp\cap \T{\p,\q}}$ 
  satisfying
  \begin{align*}
    (\c, \n_0) = (\dot{\c},\n_0) = 0,
  \end{align*}
  where $\dot{\c}$ is the derivative of $\c$ with respect to the parameter $t$.
  It is then straightforward to show that $\n(\theta, \cdot)=\ro{}{\theta}\n_0$
  is the tangent plane congruence of $\f$, hence, two rotational sphere 
  congruences span $\Lambda = \spn{\f, \n}$. Further, the principal 
  curvatures of $\f$ only depend on $t$, hence the curvature spheres 
  $s_i(\theta, \cdot) = \ro{}{\theta}\spn{\n_0 + \kappa_i \c}$ are rotational.
  
  Finally, as mentioned in the proof of \Cref{lem:ChristoffelDualLifts}, an 
  $\Omega$-surface can be characterized by the existence of lifts 
  $\s_i \in s_i$ of the curvature sphere congruences such that
  \begin{align*}
    (\s_1)_\theta=\phi_\theta \s_2, \quad (\s_2)_t = \varepsilon^2 \phi_t \s_1,
  \end{align*}
  where $\varepsilon \in \{1, i\}$. Then $\s^\pm = \s_1 \pm \varepsilon \s_2$,
  the isothermic sphere congruences, are also rotational, because 
  ${\phi_\theta=0}$. 
\end{proof}

\begin{rem}
  More generally, the curvature spheres and isothermic sphere congruences of a 
  rotational $\Omega$-surface are rotational. The above proof, written in terms
  of any rotational sphere congruences spanning $\Lambda$, can be refined to 
  show this. 
\end{rem}

%% file: Content/RotLW.tex
The goal of this section is the derivation of the differential equations 
describing rotational lW surfaces in non-flat space forms. We will restrict our 
attention to lW surfaces with real sphere complexes, where we can use 
\Cref{lem:ChristoffelDualLifts}. Our formalism will include ambient spaces of 
arbitrary (non-vanishing) constant sectional curvature and (in the case of 
hyperbolic space forms) all types of rotation. The case of parabolic rotations 
in $\H^3(\kappa)$ will be dealt with separately, because the setup for this case
involves lightlike vectors and thus degenerates to some extend. It turns out 
that the Christoffel dual lifts of rotational lW surfaces are determined by 
Jacobi elliptic differential equations. The main results of this section are 
\Cref{prop:JEF,prop:parJEF} where the saught after equations are given. These 
results will prove particularly useful in the case of cmc surfaces, investigated
in Section~\ref{sec:channelCMC}.

Let $(\Lambda, \p^\pm)$ be a rotational lW surface with real sphere 
complexes $\p^\pm$. We will describe explicit parametrizations of the isothermic
sphere congruences $s^\pm$ in terms of coordinate functions that satisfy a 
system of equations posed in \Cref{prop:JEF,prop:parJEF}. 

Let us start with the non-isotropic cases. Assume that 
the $1$-parameter subgroup of rotations $\ro{1}{}$ that fixes $\Lambda$ is 
non-parabolic. Then, we have the orthogonal splitting of the ambient space
${\R^{4,2} = \Pi \oplus_\perp \Pi_1 \oplus_\perp \Pi_2}$, where 
$\Pi = \spn{\p^+,\p^-}$, $\Pi_1$ is the rotation plane of $\rho_1$ and $\Pi_2$ 
is the orthogonal complement of $\Pi\oplus \Pi_1$. For $i=1,2$, we may choose an
orthonormal basis $\{\v_i, \e_i\}$ for $\Pi_i$ such that $\e_i$ is unit 
spacelike and $\v_i$ carries the signature of $\Pi_i$. This results in the basis
$\{\p^+, \p^-, \v_1, \e_1, \v_2, \e_2\}$ of $\R^{4,2}$ satisfying the 
multiplication table
\begin{equation}\label{eq:multiplicationTableNonIso}
  \begin{array}{c|cccccc}
    &\p^+ &\p^- &\v_1 &\e_1 &\v_2 &\e_2 \\
    \hline
    \p^+ &\kappa^+ &\varepsilon & & & & \\
    \p^- &\varepsilon &\kappa^- & & & & \\
    \v_1 & & &\kappa_1 &0 & & \\
    \e_1 & & &0 &1 & & \\
    \v_2 & & & & &\kappa_2 &0 \\
    \e_2 & & & & &0 &1
  \end{array}
\end{equation}
Given this basis, a the parametrization of a lift of a ($\ro{1}{}$-)rotational 
sphere congruence is given by
\begin{equation*}
  \s(\theta, t) = \ro{1}{\theta}(a(t), b(t), r(t), 0, h(t), k(t))^\tra,
\end{equation*}
with suitable functions $a, b, r, h, k$ parametrizing a lift of 
the profile curve $\c$ of $s$, see Subsection~\ref{sec:rotationalSurfaces}. We 
call $r, h, k$ the \emph{coordinate functions} of the surface (or of the planar 
profile curve). Further we define the constant 
$\Delta = \kappa^+\kappa^- - \varepsilon^2$.

Let $s^\pm$ be the isothermic sphere congruences enveloped by a rotational
linear Weingarten surface $(\Lambda, \p^\pm)$. According to 
\Cref{lem:rotationalLegendreLift}, $s^\pm$ and their Christoffel dual lifts 
are rotational as well. We fix Christoffel dual lifts $\s^\pm \in \Gamma s^\pm$ 
such that\footnote{Note that $\Delta \neq 0$ in the non-istropic space forms.}
\begin{equation*}
  (\s^\pm, \p^\pm) = 0, \quad \textrm{and } (\s^\pm, \p^\mp) = \Delta.
\end{equation*}
This yields the following parametrizations for $\s^\pm$: 
\begin{align}\label{eq:isoSphereParam}
  \begin{split}
    \s^+(\theta, t) &= 
      \ro{1}{\theta}(-\varepsilon, \kappa^+, r^+(t), 0, h^+(t), k^+(t))^\tra\\
    \s^-(\theta, t) &= 
      \ro{1}{\theta}(\kappa^-, -\varepsilon, r^-(t), 0, h^-(t), k^-(t))^\tra,
  \end{split}
\end{align}
for suitable functions $r^\pm, h^\pm$ and $k^\pm$. We will focus on deriving
an equation for $r^+$ and omit the decorations for the coordinate functions of 
$\s^+$ from now on. 

We define \emph{polar coordinates}: for suitable functions $\psi, d$ (of $t$) we
have
\begin{align*}
  h\v_2 + k\e_2 = d~\ro{2}{\psi}\v_2,
\end{align*}
where $\ro{2}{}$ is a $1$-parameter subgroup of rotations in $\Pi_2$. 

\begin{rem}\label{rem:polarCoordinates}
Note that this form assumes that the $\Pi_2$-part of $\s^+$ has the same causal
character as $\v^2$. We will see in Section~\ref{sec:channelCMC} that this
is the case for cmc surfaces. We call polar coordinates of this form 
\emph{causal polar coordinates} and those of the form
$h\v_2 + k\e_2 = d~\ro{2}{\psi}\e_2$ \emph{spacelike polar coordinates}.
\end{rem}

Thus, the 
parametrization of $\s^+$ in \eqref{eq:isoSphereParam} takes the form
\begin{align}\label{eq:isoSphereParamPolar}
  \s^+(\theta, t) = \ro{1}{\theta}\ro{2}{\psi(t)}(-\varepsilon, \kappa^+,
    r(t), 0, d(t), 0)^\tra
\end{align}
At this point, note that we fixed the parametrization up to a choice of
speed 
\begin{align}
  v^2 = \kappa_1 (r')^2 + \kappa_2 (h')^2 + (k')^2
\end{align}
of the profile curve.

With this notation, we pose the following proposition.

\begin{prop}\label{prop:JEF}
  The sphere congruence $\s^+$ can be parametrized such that its coordinate
  functions in \eqref{eq:isoSphereParamPolar} satisfy the following system of 
  equations:
  \begin{align}
    (r')^2 &=\frac{\kappa^-}{\Delta^2}
     (r^2 - C_1)(r^2 - C_2) \label{eq:rJEF} \\
    \psi'  &= \sqrt{\frac{\kappa_1\kappa_2}{\Delta}}
    \frac{\varepsilon r^2+\kappa^+C}{\kappa^+\Delta + \kappa_1 r^2} \label{eq:psiJEF} \\
    d^2 &= -\frac{\kappa^+\Delta + \kappa_1 r^2}{\kappa_2}, \label{eq:dJEF}
  \end{align}
  where $C_1, C_2$ are the (real) roots of the polynomial
  \begin{align}\label{eq:polynomJEF}
    p_C(x) = \kappa_1\kappa^- x^2 + (\Delta^2 + 2\kappa_1 \varepsilon C)x
      +\kappa_1\kappa^+ C^2 
  \end{align}
\end{prop}

\begin{rem}
  Note that $\sgn\Delta$ is negative if $\Pi$ is Lorentzian. Thus, 
  $\tfrac{\kappa_1\kappa_2}{\Delta}$ is always positive and \eqref{eq:psiJEF}
  allows for real solutions. Further, the specific form of \eqref{eq:rJEF} 
  assumes $\kappa^-\neq 0$. This case shall be considered separately in
  \Cref{bsp:Bryant}.
\end{rem}

\begin{proof}
  First, we note that the second equation in \eqref{eq:ChristoffelDual}
  yields\footnote{
  Note that $|\s^+_{\theta}|^2$ and $|\s^+_t|^2$ are functions of $t$ only. 
  According to Remark \ref{rem:ChristoffelDual}, the function $U$ (now a
  function of $\theta$ only) is thus constant.
  }
  \begin{align*}
    \s^-_\theta = C \frac{\s^+_\theta}{|\s^+_\theta|^2} \Rightarrow r^- = \frac{C}{r}.
  \end{align*}
  This further implies the contact condition $(\s^-,\s^+_\theta) =0$. 
  Before we introduce polar coordinates, we use that the other contact 
  conditions 
  \begin{align*}
    (\s^-, \s^+) = (\s^-, \s^+_t) = 0,
  \end{align*}
  become two equations for the two functions $h^-, k^-$, which
  we therefore express in terms of $r, h$ and $k$:
  \begin{align*}
    \begin{pmatrix}
        h^- \cr k^-
      \end{pmatrix} &= 
      \frac{1}{\kappa_2(hk'-kh')} \begin{pmatrix}
        k' &-k \cr -\kappa_2 h' &\kappa_2 h
      \end{pmatrix} \begin{pmatrix}
        -\kappa_1 C+\Delta\varepsilon  \cr
        -\frac{\kappa_1 C r'}{r}
    \end{pmatrix},
  \end{align*}
  where we assume $hk'-kh'\neq 0$. A lengthy (but uneventful) computation 
  now yields, because $\s^-$ is lightlike, 
  \begin{align*}
    0 &= \kappa^-\Delta + \kappa_1 \frac{C^2}{r^2} + 
    \frac{1}{\kappa_2(h'k - k'h)^2}\left\{
    v^2(\kappa_1 C-\Delta\varepsilon)^2 
   -r'^2\kappa_1\varepsilon^2\Delta^2
    -\kappa^+\Delta\left(\frac{\kappa_1 C r'}{r}\right)^2\right\}.
  \end{align*}
  At this point, we introduce polar coordinates for $\s^+$ to express $h$ and 
  $k$ and arrive at
  \begin{align}\label{eq:proofJEFPreliminary}
    0 =\kappa_1\Delta^3(r')^2 
      + v^2 (\kappa_1 C - \varepsilon \Delta)^2
      +\kappa_2 d^2 (\kappa^- \Delta v^2 +\kappa_1C^2\frac{v^2}{r^2}).
  \end{align}
  Note that, because $\s^+$ is lightlike,
  \begin{align*}
    \kappa_2 d^2 = -\kappa^+ \Delta - \kappa_1 r^2,
  \end{align*}
  hence we already obtained \eqref{eq:dJEF} and \eqref{eq:proofJEFPreliminary} 
  only depends on $r$ and $v$. We fix the speed of the profile curve via 
  $v^2 = r^2$ (thus introducing isothermic coordinates for $\s^+$) and finally 
  obtain
  \begin{align}\label{eq:rExpandedForm}
      (r')^2
    &=\frac{1}{\kappa_1 \Delta^2}\left\{
      \kappa_1\kappa^- r^4 + (\Delta^2 + 2\kappa_1 \varepsilon C)r^2
      +\kappa_1\kappa^+ C^2
    \right\} = \frac{1}{\kappa_1 \Delta^2} p_C(r^2).
  \end{align}
  Using the roots $C_{1,2}$ of the polynomial $p_C$ we arrive at \eqref{eq:rJEF}.
  
  Finally, note that in polar coordinates, 
  \begin{align*}
    r^2 = v^2 =
      \frac{-\kappa_1\kappa^+\Delta(r')^2+\kappa_2(\psi'd^2)^2}{\kappa_2 d^2},
  \end{align*}
  hence, \eqref{eq:rExpandedForm} and \eqref{eq:dJEF} imply \eqref{eq:psiJEF}.
\end{proof}

\begin{rems}[Solution of the equations]\label{rem:SolutionOfEquations}
  \leavevmode
  \begin{enumerate}[(a)]
    \item Equation \eqref{eq:rJEF} (as \eqref{eq:parRJEF} 
      later) is solved by Jacobi elliptic functions. 
      However, the solutions are only real if $p_C$ has real roots. Thus, we 
      first need to investigate the roots of $p_C$ whose signs will determine 
      the form of the solution to \eqref{eq:rJEF}.
    \item Note that $\psi'$ is only determined up to sign. This stems from the 
      fact that $\psi\mapsto-\psi$ amounts to a rigid motion of the surface 
      given in polar coordinates. 
  \end{enumerate}
\end{rems}

Next, we consider parabolic rotational surfaces in hyperbolic space. In this 
case $\Pi_1$ has signature $(+0)$
and is spanned by a pair of orthogonal vectors $\e_1, \v$ where $\e_1$ is unit
spacelike and $\v$ is lightlike. Chose a lightlike vector $\o$, perpendicular
to $\Pi$, such that $(\v, \o) = -1$ and a spacelike unit vector $\e_2$ to
complete a pseudo-orthonormal basis $\{\p,\q, \v, \e_1, \o, \e_2\}$ of 
$\R^{4,2}$. The multiplication table then reads as in 
\eqref{eq:multiplicationTableNonIso} except
\begin{align*}
  \begin{array}{c|cc}
     &\v &\o \\
    \hline
    \v &0 &-1 \\
    \o &-1 &0
  \end{array}.
\end{align*}
The action of the parabolic rotation $\ro{i}{}$ ($i=1,2$) is then given by
\begin{align*}
  \ro{i}{\theta}: (\v, \e_i, \o) \mapsto
    (\v_i, \e_i + \theta \v, \o + \theta \e_i + \tfrac{\theta^2}{2}\o).
\end{align*}
This implies the following parametrizations for Christoffel 
dual lifts $\s^\pm$ of isothermic sphere congruences enveloped by a parabolic
rotational linear Weingarten surface $(\Lambda, \p^\pm)$
\begin{align}\label{eq:parIsoSphereParam}
  \s^+(\theta, t) &= 
  \ro{1}{\theta}\left(-\varepsilon, \kappa^+, h(t),0,r(t),k(t)\right)^\tra \\
  \s^-(\theta, t) &=
  \ro{1}{\theta}\left(\kappa^-, -\varepsilon, h^-(t), 0, r^-(t), k^-(t)\right)^\tra.
\end{align}
These lifts are again uniquely determined by requiring 
$(\s^\pm, \p^\mp) = \Delta$, however via reparametrization we can still choose 
the speed 
\begin{align*}
  v^2 = (k')^2 -2r'h',
\end{align*}
of the profile curve of $\s^+$. 

Polar coordinates, in this case, are obtained via rotations $\ro{2}{}$ 
in the plane $\spn{\v, \e_2}$ to write
\begin{align*}
  k(t)\e_2 + h(t)\o = \ro{2}{\psi(t)}\o + d(t) \v.
\end{align*}
This yields the polar coordinate parametrization
\begin{align}\label{eq:parIsoSphereParamPolar}
  \s^+(\theta, t) = \ro{1}{\theta} \ro{2}{\psi(t)} \left(
    -\varepsilon, \kappa^+, d(t), 0, r(t), 0)
  \right),
\end{align}
for the Christoffel dual lift $\s^+ \in \Gamma s^+$.

\begin{prop}\label{prop:parJEF}
  The sphere congruence $\s^+$ can be parametrized such that its coordinate
  functions in \eqref{eq:parIsoSphereParamPolar} satisfy the following system of
  equations:
  \begin{align}
    (r')^2 &= \frac{\kappa^-}{\Delta^2}(r^2 - C_1)(r^2-C_2),\label{eq:parRJEF} \\
    \psi'  &= \sqrt{-\tfrac{1}{\Delta}}\frac{\varepsilon r^2 + \kappa^+ C}{r^2}\label{eq:parPsiJEF} \\
    d &=  \frac{\kappa^+ \Delta}{2r}, \label{eq:parDJEF}
  \end{align}
   where $C_1, C_2$ are the (real) roots of the polynomial
  \begin{align}\label{eq:parPolynomJEF}
    p_C(x) = \kappa^- x^2 + 2C\varepsilon x +\kappa^+ C^2.
  \end{align}
\end{prop}

The proof of this proposition follows the same strategy as the one for 
\Cref{prop:JEF}. 

\begin{bsp}\label{bsp:Bryant}
  The form we chose for \eqref{eq:rJEF} and \eqref{eq:parRJEF} assumes 
  that the polynomial $p_C$ is indeed quadratic, hence $\kappa^- \neq 0$. 
  We will consider the case of $\p^-$ being lightlike here: this 
  corresponds to the class of either minimal surfaces in $\R^3$ or 
  Bryant-type surfaces in $\H^3$ \cite{galvez2004}. In the case of Bryant type 
  surfaces, the first part of the proof of \Cref{prop:JEF} applies, but 
  \eqref{eq:rExpandedForm} becomes
  \begin{align*}
    (r')^2 =
      (\frac{1}{\kappa_1} +2 \frac{\varepsilon C}{\Delta^2})r^2
      +\frac{\kappa^+}{\Delta^2} C^2.
  \end{align*}
  The solution to this depends on the signs of the constants $\kappa^+$ and 
  $A:=\tfrac{1}{\kappa_1} +2 \tfrac{\varepsilon C}{\Delta^2}$: 
  \begin{align*}
    \begin{array}{ll}
      \kappa^+ < -,\, A<0:
        &\textrm{no solution} \\
      \kappa^+ < 0,\, A>0: 
        &r(t) = \sqrt{\tfrac{-\kappa^+ C^2}{A \Delta^2}}\,\cosh(\sqrt{A}\,t) \\
      \kappa^+ > 0,\, A<0:
        &r(t) = \sqrt{\frac{\kappa^+ C^2}{-A\Delta^2}}\,\sin(\sqrt{-A}\,t) \\
      \kappa^+ > 0,\, A>0:
        &r(t) = \sqrt{\tfrac{\kappa^+ C^2}{A\Delta^2}}\,\sinh(\sqrt{A}\,t).
    \end{array}
  \end{align*}
  For parabolic rotations, the solution space is the same, but 
  $A:=2\tfrac{\varepsilon C}{\Delta^2}$. The function $\psi$ is found
  via simple integration of these trigonometric/hyperbolic functions.
\end{bsp}

%% file: Content/RotCMC.tex
In this section, we apply the results of Section~\ref{sec:rotLW} to
the specific case of rotational cmc surfaces. We will split our
investigation into two parts, considering surfaces with $H^2+\kappa>0$
(which we call \emph{Delaunay-type surfaces}) and surfaces with $H^2+\kappa<0$
(called \emph{sub-horospherical surfaces}) separately. Due to the similarity
of the first class to cmc surfaces in $\R^3$, which lacks for the latter class,
this separation comes naturally and will produce differences
in the solution (of \eqref{eq:rJEF}) and geometric properties of the results. 
Further, this section will include application of the explicit parametrizations 
we obtain to recover computational proves of some well-known facts about 
rotational cmc surfaces (see \Cref{thm:cmcTori}).

From now on, we assume that for cmc surfaces $H\geq 0$. Given any cmc surface, 
this can be achieved via a change of orientation of the Gauss map.

\subsection{Delaunay-type surfaces}\label{sec:DelaunayCase}
\input{Content/ChannelCMC/Delaunay}

\subsection{Sub-horospherical surfaces in
 \texorpdfstring{$\H^3$}{H3}}\label{sec:smallH}
\input{Content/ChannelCMC/SmallH}

%% file: Content/ChannelCMC/Delaunay.tex
In this subsection, we apply the equations derived in 
Subsection~\ref{sec:rotLW} to the case of rotational Delaunay-type surfaces
in $\S^3$ and $\H^3$. Delaunay-type surfaces, that is, rotational cmc surfaces 
with $H^2+\kappa\geq 0$, have a rich history and have been described from 
multiple viewpoints, see for instance \cite{bryant1987,dorfmeister2000}.

Furthermore, these surfaces are parallel to cGc surfaces, hence one can obtain
parametrizations in terms of Jacobi elliptic functions via parallel 
transformation from the parametrizations given for rotational cGc surfaces in 
\cite{hertrich-jeromin2023}. However, we shall use the equations derived in
Subsection~\ref{sec:rotLW}, because this approach is
\begin{inparaenum}[1)]
  \item more direct, 
  \item particularly useful in the cmc case and
  \item also applicable for surfaces with $H^2+\kappa\leq 0$ (which only exist 
    in $\H^3$).
\end{inparaenum}

As we have mentioned in \Cref{bsp:conservedQuantities}, the sphere complexes of 
a cmc surface $(\Lambda, \p^\pm)$ in a space form $\Q{\p,\q}$ take the form
\begin{align}\label{eq:cmcConservedQuantities}
  \p^+=\p, \quad \p^-=\q - H \p
\end{align}
and Christoffel dual lifts of the isothermic sphere congruences $s^\pm$ are 
given as combinations of the space form lifts $\f, \n$ of the point sphere 
envelope and the tangent plane congruence via
\begin{align}\label{eq:cmcIsothermicSphereCong}
  \s^+ = \f, \quad \s^- = \n + H \f,
\end{align}
up to constant rescalings.

\begin{rem}
  The fact that $s^+$ represents the surface itself implies that we have to
  use causal polar coordaintes (see Remark \ref{rem:polarCoordinates}): whenever
  $\Pi_2$ has signature $(++)$ this is obvious. For signature $(+-)$ (which 
  only occurs in hyperbolic space forms), the $\Pi_2$-part of $\s^+$ has
  to be timelike, or $\s^+$ would be spacelike. 
\end{rem}

To apply \Cref{prop:JEF} (or \Cref{prop:parJEF} in the parabolic rotation case),
we note 
\begin{align}\label{eq:quantitiesValueDelaunay}
  \kappa^+ = -1,~ \kappa^- = -(H^2 +\kappa), ~\varepsilon = H~\Rightarrow ~
  \Delta = \kappa.
\end{align}
Since $(\s^+, \p^-) = \kappa$, the Christoffel dual lift described in
\Cref{prop:JEF} coincides with the space form lift $\f$ for $\kappa\in\{\pm 1\}$
(up to sign). This explains why the equations derived in 
Subsection~\ref{sec:rotLW} are particularly useful in the case of cmc surfaces. 

\begin{rem}
  Because $\s^+$ parametrizes a surface in a space form, we will utilize 
  polar coordinates that are not spacelike polar coordinates (see 
  Remark~\ref{rem:polarCoordinates}). 
\end{rem}

\begin{bsp}\label{bsp:BryantCMC}
  The surfaces considered in \Cref{bsp:Bryant} contain the class of cmc
  surfaces in hyperbolic space with ${\kappa^- = H^2+\kappa = 0}$. 
  These surfaces, dubbed \emph{horospherical} in 
  \cite{hertrich-jeromin2000}, are thus parametrized by the coordinate
  functions
  \begin{align*}
      &r(t) = \sqrt{\tfrac{C^2}{A}}\,\cosh(\sqrt{A}\,t), 
      &\psi(t) = \sqrt{-\frac{\kappa_1}{\kappa_2}}\frac{1}{\kappa_1}
      \left[
        t + \frac{\kappa_1 C^2 - 1}{\sqrt{A+\kappa_1 C^2}}
        \operatorname{artanh}\left(
        \frac{\operatorname{tanh}(\sqrt{A}~t)}{\sqrt{1+\tfrac{\kappa_1 C^2}{A}}}
        \right)
      \right]
  \end{align*}
  where $A=\tfrac{1}{\kappa_1} +2 C$ needs to be positive, 
  hence $C > -\tfrac{1}{2\kappa_1}$ (%
  $A=2C$ and $C>0$ in the parabolic case). This gives (upon integration of
  \eqref{eq:parPsiJEF} in the parabolic case), a complete classification of 
  rotational horospherical surfaces including, for instance, catenoid cousins 
  (compare \cite{bryant1987,umehara1993}).
\end{bsp}

We start with the case of non-parabolic surfaces and investigate the polynomial
\eqref{eq:polynomJEF}: for $H^2+\kappa >0$, 
\begin{align*}
  p_C(x)= -\kappa_1 (H^2 + \kappa) x^2 + (\kappa^2 + 2\kappa_1 H C)x -
      \kappa_1C^2,
\end{align*}
has real roots $C_1, C_2$ if and only if
\begin{align*}
  D(C):=(\kappa^2 + 2\kappa_1 HC)^2 -4\kappa_1^2 C^2 (H^2 +\kappa) >0. 
\end{align*}
This poses a condition on the constant $C$, which we shall investigate case 
by case. However, in all cases we have
\begin{align*}
  C_1 C_2 = \frac{C^2}{H^2 + \kappa} \geq 0,
\end{align*}
with $C=0$ yielding an equal sign. This implies that $C_1$ and $C_2$ have the
same sign. However, under the assumption $C_1, C_2 <0$, \eqref{eq:rJEF} becomes
\begin{align*}
  (r')^2 = -\frac{H^2 + \kappa}{\kappa^2}\left(r^2 + (-C_1)\right)
    \left(r^2 + (-C_2)\right),
\end{align*}
which allows for no real solutions. Thus, we further restrict $C$ to values 
where $C_1, C_2 \geq 0$. From now on, we name the roots such that 
$C_1 > C_2 \geq0$. 

Note that 
\begin{align*}
  D(x)=0 ~ \Leftrightarrow ~ x^\pm=\frac{\kappa}{2\kappa_1}\left(
    H\pm\sqrt{H^2 + \kappa}
  \right),
\end{align*}
which determines the conditions on $C$: 
\begin{itemize}
	\item[$\kappa>0, \kappa_1 >0$] In this case, $C \in [x^-, x^+]$ guarantees 
    real roots of $p_C(x)$. Note that $C_1 > 0$ at ${C=0\in [x^-, x^+]}$, hence 
    for all feasable values, $C_1 > C_2 \geq 0$ and at the bounds $C_1 = C_2$. 
  \item[$\kappa<0, \kappa_1 > 0$] In this case, 
    $C\in (-\infty, x^-] \cup [x^+, \infty)$ guarantees real roots. However, 
    $0 \in [x^+, \infty)$, which we know to be a valid value. Indeed, for 
    $C \in (-\infty, x^-]$, we have $C_1, C_2 <0$. We conclude that we restrict
    $C$ to $C\in[x^+, \infty)$ and have $C_1 > C_2 \geq 0$ and $C_1 = C_2$ 
    at $C=x^+$.
  \item[$\kappa<0, \kappa_1 <0$] In the case of hyperbolic rotations in $\H^3$,
    finally, $C \in (-\infty, x^+] \cup [x^-, \infty)$. In this case, however, 
    $C=0$, is not a feasible value. Thus, we conclude $C \in [x^-, \infty)$, 
    which yields $C_1 > C_2 \geq 0$ and $C_1 = C_2$ at $C=x^-$.  
\end{itemize}

With these preparations, we prove in the following proposition that the solution
$r$ of \eqref{eq:rJEF} can be expressed as a Jacobi elliptic function and the 
solution $\psi$ of \eqref{eq:psiJEF} can be expressed as an elliptic integral: 
$\Pii{k}{p}{s}$ denotes the incomplete elliptic integral of 
the third kind\footnote{Defining the incomplete integral of third 
kind as
\begin{align*}
\Pi(k;s,p)=
\int_0^{s}\tfrac{1}{1-k\sin^2(u)}\tfrac{du}{\sqrt{1-p^2\sin^2(u)}},
\end{align*}
as is often done, we obtain the relationship
\begin{align*}
	\Pii{p}{k}{s} = \Pi(k; \jac{am}{p}(s), p).
\end{align*}
} with modulus $p$ and parameter $k$ as defined in 
\cite[Sec 17.2]{abramowitz1972}, that is, 
\begin{align*}
	\Pii{p}{k}{s} = \int_0^{s}\tfrac{du}{1-k \jac{sn}{p}^2(u)}.
\end{align*}

\begin{prop}\label{prop:rOfDelaunay}
  Assume $H^2 + \kappa>0$ and let $C$ be chosen such that the polynomial $p_C$
  has two real roots $C_1, C_2$ satisfying $C_1> C_2 \geq 0$. Then, the solutions
  of \eqref{eq:rJEF} and \eqref{eq:psiJEF} are given by
  \begin{align}\label{eq:JEFSolutionDelaunay}
    \begin{split}
      r(t) &= \sqrt{C_1} \jac{dn}{p}\left(
      \FAC\,t
      \right), \textrm{ with } p^2 = \frac{C_1-C_2}{C_1} \in [0,1], \\
      \psi(t) &= \sqrt{\frac{\kappa_1\kappa_2}{\kappa}}\frac{1}{\kappa_1}
      \left\{Ht - \frac{\kappa_1 C - \kappa H}{\kappa_1 C_1 - \kappa}
      \frac{1}{\FAC}
      \Pii{p}{\frac{\kappa_1 (C_1-C_2)}{\kappa_1 C_1 - \kappa}}{\FAC\,t}\right\},
    \end{split}
  \end{align}
  where $\FAC = \sqrt{\tfrac{H^2 + \kappa}{\kappa^2}C_1}$.
\end{prop}

\begin{proof}
  Since $C_1 > 0$ we define a (real-valued) function $y$ satisfying
  $r=\sqrt{C_1}y$. We then have
  \begin{align*}
    (r')^2 &= -\frac{H^2 + \kappa}{\kappa^2}(r^2 - C_1)(r^2 - C_2) \\
    \Leftrightarrow 
    (y')^2 &= \frac{H^2 + \kappa}{\kappa^2}\,C_1(1-y^2)(-\frac{C_2}{C_1} + y^2),
  \end{align*}
  the real solution of which is
  \begin{align*}
    y(t) = \jac{dn}{p}\left(
      \sqrt{\frac{H^2 + \kappa}{\kappa^2}\,C_1}\, t
    \right), \textrm{ with } p^2 = \frac{C_1-C_2}{C_1}.
  \end{align*}
  Note that $p\in [0,1]$ because $C_1 > C_2$. This proves to proclaimed form of 
  $r$. 
  
  For $\psi$, we put the expression for $r$ into \eqref{eq:psiJEF}. Using the 
  Jacobi identity $\jac{dn}{p}^2 + p^2 \jac{sn}{p}^2 = 1$, this yields
  \begin{align*}
    \psi'(t)&= \sqrt{\frac{\kappa_1\kappa_2}{\kappa}}
      \frac{H C_1 \jac{dn}{p}^2\left(
      \sqrt{\tfrac{H^2 + \kappa}{\kappa^2}\,C_1}\, t
    \right) - C}{-\kappa + \kappa_1 C_1 \jac{dn}{p}^2\left(
      \sqrt{\tfrac{H^2 + \kappa}{\kappa^2}\,C_1}\, t
    \right)} \\
      &=\sqrt{\frac{\kappa_1\kappa_2}{\kappa}}\frac{1}{\kappa_1}\left\{H -
      \frac{\kappa_1 C - \kappa H}{\kappa_1 C_1 - \kappa}
      \frac{1}
      {1 -\frac{\kappa_1 (C_1-C_2)}{\kappa_1 C_1 - \kappa} \jac{sn}{p}^2\left(
      \sqrt{\tfrac{H^2 + \kappa}{\kappa^2}\,C_1}\, t
    \right)} \right\}
  \end{align*}
  By integration we obtain the desired solution $\psi$. 
\end{proof}

\begin{rem}
  The given form for $\psi$ is invalid if $\kappa_1 C_1 - \kappa = 0$. However, 
  this happens if and only if $\kappa_1C - \kappa H = 0$ in which case 
  $\psi(t) = Ht$ as one quickly deduces from \eqref{eq:psiJEF} (up to sign 
  ambiguity). We will not emphasize this special case from now on.
\end{rem}

We now apply this result to give explicit parametrizations of Delaunay-type 
surfaces in $\S^3$ and $\H^3$ in the following two theorems. These directly 
follow from \Cref{prop:rOfDelaunay} under the right choice of constants (which
we will highlight).

We regard $\S^3$ as the unit sphere in $\R^4$, which is equipped with an 
orthonormal basis. This corresponds to $\kappa = \kappa_1 = \kappa_2 = 1$ in
\Cref{prop:rOfDelaunay}.

\begin{thm}[Classification of Delaunay-type surfaces in 
  $\S^3$]\label{thm:DelaunayS3}
  Every rotational constant mean curvature $H$ surface in 
  $\S^3 \subset \R^4$ is given by
  \begin{align*}
    (\theta, t) \mapsto \left(
      r(t) \cos \theta, r(t) \sin \theta, 
      \sqrt{1-r^2(t)} \cos \psi(t), \sqrt{1-r^2(t)} \sin \psi(t),
    \right)
  \end{align*}
  with 
  \begin{align*}
    r(t) &= \sqrt{C_1} \jac{dn}{p}\left(
      \FAC\,t
      \right), \textrm{ with } p^2 = \frac{C_1-C_2}{C_1} \in [0,1], \\
    \psi(t) &= Ht - \frac{C - H}{C_1 - 1}
      \frac{1}{\FAC}
      \Pii{p}{\frac{C_1-C_2}{C_1 - 1}}{\FAC\,t},
  \end{align*}
  where $\FAC = \sqrt{(H^2 + 1)C_1}$ and $C_1>C_2$ are
  the non-negative roots of 
  \begin{align*}
      p_C(x)= (H^2 + 1) x^2 - (1 + 2 H C)x + C^2.
  \end{align*}
\end{thm}

For the following theorem, we view $\H^3$ in the Minkowski model, that is, as a 
subset in $\R^{3,1}$, which is equipped with an orthonormal basis 
$(\e_0, \e_1,\e_2, \e_3)$, where $\e_0$ is unit timelike. This amounts to 
$\kappa=-1$, $\kappa_1 = - \kappa_2$  and $\kappa_1 = 1$ for elliptic and
$\kappa_1 = -1$ for hyperbolic rotations.

\begin{thm}[Classification of non-parabolic Delaunay-type 
  surfaces in $\H^3$]\label{thm:DelaunayH3}
  Every non-parabolic rotational constant mean curvature surface with
  $H^2>1$ in $\H^3 \subset \R^{3,1}$ is given by one of the following 
  parametrizations: 
  \begin{itemize}
    \item[elliptic rotations: $\kappa_1=1$] \begin{align*}
        (\theta, t)\mapsto \left(
          \sqrt{1 + r^2(t)} \cosh \psi(t), 
          \sqrt{1 + r^2(t)} \sinh \psi(t),
          r(t) \cos\theta, 
          r(t) \sin\theta
        \right),
      \end{align*}
    \item[hyperbolic rotations: $\kappa_1=-1$] \begin{align*}
        (\theta, t) \mapsto \left(
          r(t) \cosh\theta, 
          r(t) \sinh\theta, 
          \sqrt{r^2(t)-1} \cos\psi(t), 
          \sqrt{r^2(t)-1} \sin\psi(t)
        \right)
      \end{align*}
  \end{itemize}
  with 
    \begin{align*}
      r(t) &= \sqrt{C_1} \jac{dn}{p}\left(
        \FAC\,t
      \right), \textrm{ with } p^2 = \frac{C_1-C_2}{C_1} \in [0,1], \\
      \psi(t) &=Ht - \frac{C + \kappa_1 H}{C_1 + \kappa_1}
        \frac{1}{\FAC}
        \Pii{p}{\frac{C_1-C_2}{C_1 + \kappa_1}}{\FAC\,t},
    \end{align*}
    where $\FAC = \sqrt{(H^2 -1)C_1}$ and $C_1 > C_2$ are the non-negative
    roots of 
    \begin{align*}
      p_C(x) = (H^2 - 1) x^2 - (\kappa_1 + 2 H C)x + C^2.
    \end{align*}
\end{thm}

In order to finish the classification, we consider parabolic Delaunay-type 
surfaces in $\H^3$. Again, we regard $\H^3$ as subset of $\R^{3,1}$, however we 
equip $\R^{3,1}$ with a pseudo-orthonormal basis $(\v, \e_1, \o, \e_2)$. This 
corresponds to $\kappa=-1$ in \Cref{prop:parJEF}. Solving the differential 
equations therein leads to \Cref{thm:parDelaunayH3}. First, however, we discuss 
the polynomial $p_C$ for this case.

The roots of $p_C$ take the particularly simple form
\begin{align*}
  C^\pm = \frac{C}{H\pm 1}.
\end{align*}
For negative $C$, these are both negative, hence, no real solution of 
\eqref{eq:parRJEF} exists. For positive $C$, we have $C^- > C^+ > 0$\footnote{
For $C=0$, $p_C$ degenerates and \eqref{eq:parRJEF} does not have real 
solutions.}. 
The solutions of \eqref{eq:parRJEF} and \eqref{eq:parPsiJEF} are now obtained
in a way similar to the proof of \Cref{prop:rOfDelaunay}. Thus, we arrive at 
the following theorem.

\begin{thm}[Classification of parabolic Delaunay-type 
  surfaces in $\H^3$]\label{thm:parDelaunayH3}
  Every parabolic rotational constant mean curvature surface with
  $H^2>1$ in $\H^3 \subset \R^{3,1}$ is given by
  \begin{align*}
    (\theta, t)\mapsto \left(
      \frac{1+r^2(\psi^2 + \theta^2)}{2r},
      r(t)\theta, 
      r(t), 
      r(t)\psi(t)
    \right),
  \end{align*}
  and the functions
  \begin{align*}
    r(t)&= \sqrt{C^-}\jac{dn}{p}\left(
      \FAC\,t
    \right) \textrm{ with } p^2 = \frac{2}{H+1} \in (0,1)\\
    \psi(t)&= Ht - \frac{C}{C^-}\frac{1}{\FAC}
      \Pii{p}{p^2}{\FAC\,t},
  \end{align*}
  where $\FAC = \sqrt{C(H+1)}$ and $C^\pm = \tfrac{1}{H\pm 1}$.
\end{thm}

\begin{rem}\label{rem:unduloidNodoid}
As in the Euclidean case, Delaunay-type surfaces
in $\S^3$ and $\H^3$ come in two flavors. Namely, unduloids and the nodoids can 
be characterized by whether $\psi$ is monotone or oscillating. This
bifurcation appears in the parametrizations of 
\Cref{thm:DelaunayS3,thm:DelaunayH3} in the following way: 
unless $\psi'$ is constant, we have $A\psi' =Hr^2 - C$, where $A$ is a positive
function. This implies that the sign of $\psi'$ is constant if and only of 
$\tfrac{C}{H}$ is not in the image $[C_2, C_1]$ of $r^2$. We can plug in 
specific values for $C$ (namely $0, H$ and $2H$) to obtain situations where 
$\psi'$ has constant sign (unduloid case) or oscillates (nodoid case). In the
parabolic rotational case, all surfaces are nodoids. 
\end{rem}

As an application of our parametrizations, we give a simple construction of 
embedded cmc tori in $\S^3$. As was proved in \cite{andrews2015}, 
such surfaces are always rotational and can thus be constructed using 
\Cref{thm:DelaunayS3}. For given mean curvature value $H$ the question becomes
which values for the constant $C$ yield closed profile curves. This is answered
in the following theorem. 

\begin{thm}\label{thm:cmcTori}
  Given a Delaunay-type surface in $\S^3$, parametrized as in 
  \Cref{thm:DelaunayS3}, the profile is periodic if there exist 
  $n\in \N$ and $C\in \R$ feasible such that
  \begin{align}\label{eq:Closed}
    \frac{n}{\FAC}\left(F_p H-\frac{C-H}{C_1-1}\Pi_p\left(\frac{C_1-C_2}{C_1-1}%
    \right)\right) = \pi,
  \end{align}
  where $F_p$ and $\Pi_p(\frac{C_1-C_2}{C_1-1})$ denote the complete elliptic
  integrals of first and third kind respectively. Additionally, if $C<0$, the
  profile curve has no self-intersections and the resulting surface is an 
  embedded cmc torus of revolution.
\end{thm}

\begin{rem}
  A similar theorem can be stated for hyperbolic rotational Delaunay-type 
  surfaces in $\H^3$, where we can use a formula similar to \eqref{eq:Closed}
  to find complete and embedded hyperbolic rotational cmc surfaces.
\end{rem}

\begin{proof}
  The radius function given in \Cref{thm:DelaunayS3} is periodic with period 
  $\alpha= \tfrac{2F_p}{\FAC}$. Using this and 
  \begin{align*}
    \Pii{p}{k}{x+nF_p} = \Pii{p}{k}{x} + n\Pi_p(k),
  \end{align*}
  we see that \eqref{eq:Closed} is equivalent to
  \begin{align*}
    \psi(s+\alpha) = \psi(s) + 2\pi.
  \end{align*}
  Because $\cos\psi$ and $\sin\psi$ are thus $\alpha$-periodic, so is the
  profile curve of the surface. The fact that it has no self-intersections 
  for $C<0$ follows from the discussion in \Cref{rem:unduloidNodoid}.
\end{proof}

\begin{figure}%
\centering
\begin{tikzpicture}
	\node at (0,0) {\includegraphics[width=0.5\textwidth]{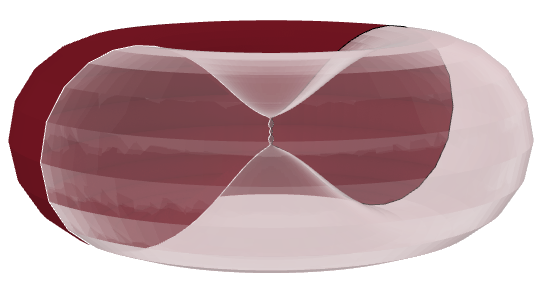}};
  \node at (8,0) {\includegraphics[width=0.25\textwidth]{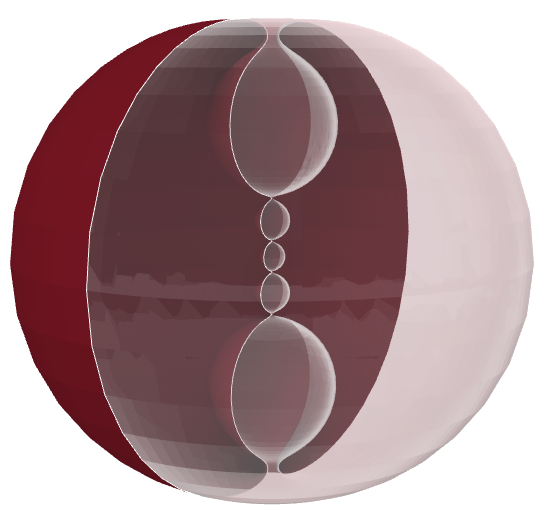}};
\end{tikzpicture}
\caption{Two cmc $H=2$ tori in $\S^3$ after stereographic projection ($n=5,6$).
Plotted and rendered using Rhino and mathematica.}%
\label{fig:cmcTori}%
\end{figure}

\begin{bsp}
  Figure~\ref{fig:cmcTori} depicts two cmc tori with $H=2$ in the sphere after
  stereographic projection into $\R^3$. These correspond to the cases $n=5, 6$ 
  (and suitable values for $C$) in \eqref{eq:Closed}.
\end{bsp}

%% file: Content/ChannelCMC/SmallH.tex
In this subsection, we apply \Cref{prop:JEF} (and \Cref{prop:parJEF}) to the 
case of surfaces with $H^2 + \kappa<0$ in $\H^3$, which we call 
\emph{sub-horospherical}. For sub-horospherical surfaces in hyperbolic space 
forms, less literature is available than for their Delaunay-type counterparts. 
This stems from the fact that they are less similar to the well-studied class of
cmc surfaces in $\R^3$, i.e., there is no Lawson correspondence, and the DPW 
method \cite{dorfmeister2000} in its original form could not be used to 
construct surfaces with $H^2 + \kappa<0$. However, the special case of 
rotational surfaces was considered in some publications catering to rotational 
surfaces of constant curvature (with any value $H$). For instance 
\cite{bobenko2019} uses a new loop group method (\cite{dorfmeister2014}) to 
construct examples of rotational cmc surfaces and \cite{dursun2020,mori1983} 
give explicit parametrizations of their profile curves.

In our setup, the significant difference
in the solutions to \eqref{eq:rJEF} arises in the form of the roots of the 
polynomial $p_C$. We will analyze these differences to obtain explicit 
parametrization of rotational sub-horospherical surfaces in terms of Jacobi 
elliptic functions. These results also apply to the well-studied class of 
\emph{minimal surfaces} ($H=0$) in hyperbolic space (see for example 
\cite{bobenko2019,kokubu1997}).

The setup is as in Subsection~\ref{sec:DelaunayCase}: we have 
\begin{align*}
  \kappa^+ = -1,~ \kappa^- = -(H^2 +\kappa), ~\varepsilon = H~\Rightarrow ~
  \Delta = \kappa<0,
\end{align*}
and \eqref{eq:polynomJEF} takes the form
\begin{align*}
  p_C(x) = -(H^2+\kappa)\kappa_1 x^2 + (\kappa^2 +2\kappa_1 H C)x -\kappa_1C^2.
\end{align*}
The main difference to Subsection~\ref{sec:DelaunayCase} lies in the following
facts: $p_C$ has real roots $C^+, C^-$ for all values of $C$ 
but they satisfy
\begin{align*}
  C^+C^- = \frac{C^2}{H^2 + \kappa} < 0,
\end{align*}
hence, with the right choice of nomenclature, $C^+ > 0 \geq C^-$ ($C^-$ vanishes
for $C=0$).

This warrants a different analysis of equation \eqref{eq:rJEF}, which we shall 
give in the proof of the following proposition.

\begin{prop}\label{prop:rOfSmallH}
  Assume $H^2 + \kappa<0$ in a hyperbolic space form and choose $C\in \R$. Then,
  the solutions of \eqref{eq:rJEF} and \eqref{eq:psiJEF} are given by
  \begin{align}\label{eq:JEFSolutionSmallH}
    \begin{split}
      r(t)&= \frac{\sqrt{C^+}}{\jac{cn}{p}\left(\FAC\,t\right)}, \\
      \psi(t)&=\sqrt{\frac{\kappa_1\kappa_2}{\kappa}}\frac{1}{\kappa}\left\{
        C t + \frac{C^+(H\kappa-C\kappa_1)}{\kappa_1 C^+ - \kappa}
        \frac{1}{\FAC}\Pii{p}{\frac{\kappa}{\kappa - \kappa_1 C^+}}{\FAC\,t}
        \right\}
    \end{split}
  \end{align}
  with $p^2 = -\frac{C^-}{C^+ - C^-} \in [0,1]$ and 
  $\FAC=\sqrt{-(H^2 + \kappa)(C^+ - C^-)}$.
\end{prop}

\begin{proof}
  The manipulation of the differential equation \eqref{eq:rJEF} follows the 
  treatment described in the proof of \Cref{prop:rOfDelaunay}: define
  $y=\tfrac{r}{\sqrt{C^+}}$, a real function that satisfies
  \begin{align*}
    (y')^2 &=\frac{H^2 + \kappa}{\kappa^2}(C^+ - C^-)(1-y^2)
               \left(-\frac{C^-}{C^+ - C^-} + \frac{C^+}{C^+ - C^-} y^2\right),
  \end{align*}
  due to \eqref{eq:rJEF}. This implies
  \begin{align*}
    y(s) = \jac{cn}{\tilde{p}}\left(i\sqrt{-\frac{H^2 + \kappa}{\kappa^2}(C^+ - C^-)}\,s\right), 
    \textrm{ with } \tilde{p}^2 = \frac{C^+}{C^+ - C^-}<0,
  \end{align*}
  hence, using Jacobi's imaginary transformation (see 
  \cite[App A]{hertrich-jeromin2023}), 
  \begin{align*}
    y(s) = \frac{1}{\jac{cn}{p}\left(\sqrt{-\frac{H^2 + \kappa}{\kappa^2}(C^+ - C^-)}\,s\right)}, 
    \textrm{ with } p^2 = -\frac{C^-}{C^+ - C^-}.
  \end{align*}
  This proves the correctness of the stated solution $r$. Clearly, $p^2$ is 
  positive but smaller than $1$ (recall that $C^- \geq 0$). We 
  can use this solution in \eqref{eq:psiJEF} to obtain  
  \begin{align*}
    \psi'(s)&=\sqrt{\frac{\kappa_1\kappa_2}{\kappa}}
        \frac{HC^+\jac{nc}{p}^2\left(\sqrt{-(H^2 + \kappa)(C^+ - C^-)}\,s\right) - C}
        {-\kappa+\kappa_1 C^+\jac{nc}{p}^2\left(\sqrt{-(H^2 + \kappa)(C^+ - C^-)}\,s\right)}\\
        &=\sqrt{\frac{\kappa_1\kappa_2}{\kappa}}\frac{1}{\kappa}\left\{
        C + \frac{C^+(H\kappa-C\kappa_1)}{\kappa_1 C^+ - \kappa}
        \frac{1}{1 - \frac{\kappa}{\kappa - \kappa_1 C^+} 
        \jac{sn}{p}^2\left(\sqrt{-(H^2 + \kappa)(C^+ - C^-)}\,s\right)}
        \right\}.
  \end{align*}
  Integration leads to the claimed form of $\psi$.
\end{proof}

As application of \Cref{prop:rOfSmallH} we give the classification of 
non-parabolic sub-horospherical surfaces in $\H^3$. We arrive at this by
fixing $\kappa=-1$ and an orthonormal basis $\{\v_1, \e_1, \v_2, \e_2\}$
of $\R^{3,1} \cong \spn{\p^\pm}^\perp$ such that 
$\kappa_1 = -\kappa_2 \in \{\pm 1\}$. 

\begin{thm}[Classification of non-parabolic sub-horospherical surfaces
  in $\H^3$]\label{thm:smallH}
  Every non-parabolic rotational constant mean curvature $H$ surface with 
  $H^2<1$ in $\H^3$ is given by one of the following parametrizations:
    \begin{itemize}
    \item[elliptic rotations: $\kappa_1=1$] \begin{align*}
        (\theta, t)\mapsto \left(
          \sqrt{1 + r^2(t)} \cosh \psi(t), 
          \sqrt{1 + r^2(t)} \sinh \psi(t),
          r(t) \cos\theta, 
          r(t) \sin\theta
        \right),
      \end{align*}
    \item[hyperbolic rotations: $\kappa_1=-1$] \begin{align*}
        (\theta, t) \mapsto \left(
          r(t) \cosh\theta, 
          r(t) \sinh\theta, 
          \sqrt{r^2(t)-1} \cos\psi(t), 
          \sqrt{r^2(t)-1} \sin\psi(t)
        \right)
      \end{align*}
  \end{itemize}
  with 
    \begin{align*}
      r(t) &=\frac{\sqrt{C^+}}{\jac{cn}{p}\left(\FAC\,t\right)}, 
        \textrm{ with } p^2 = -\frac{C^-}{C^+ - C^-} \in [0,1], \\
      \psi(t) &= 
      C t - \frac{C^+(C + \kappa_1 H)}{C^+ + \kappa_1}
      \frac{1}{\FAC}\Pii{p}{\frac{1}{\kappa_1 C^+ + 1}}{\FAC\,t},
    \end{align*}
    where $\FAC = \sqrt{-(H^2 - 1)(C^+ - C^-)}$ and $C_1 > 0 \geq C_2$ 
    are the (real) roots of 
    \begin{align*}
      p_C(x) = (H^2 - 1) x^2 - (\kappa_1 + 2 H C)x + C^2.
    \end{align*}
\end{thm}

We are left with the case of parabolic rotations: as in 
Subsection~\ref{sec:DelaunayCase}, equip 
$\R^{3,1}$ with a pseudo-orthonormal basis $(\v, \e_1, \o, \e_2)$. This 
corresponds to $\kappa=-1$ in \Cref{prop:parJEF}. The roots of the polynomial 
\eqref{eq:parPolynomJEF} are
\begin{align*}
  C_{1,2} = \frac{C}{H\pm 1}.
\end{align*}
As $H^2<1$, for all values of $C$, one of these roots is positive and one negative 
(excluding the degenerate case $C=0$). Let us denote the positive root by $C^+$ 
and the negative one by $C^-$. We then see, as in the proof of 
\Cref{thm:parDelaunayH3}, that
\begin{align*}
  r(t) = \frac{\sqrt{C^+}}{\jac{cn}{p}(\FAC\,t)},
\end{align*}
where $\FAC = \sqrt{|2C|}$ and  
\begin{align*}
  p^2 = \left\{\begin{array}{ll}
      \frac{H+1}{2} &\textrm{for } C>0 \\
      \frac{-H+1}{2} &\textrm{for } C<0.
    \end{array}\right.
\end{align*}
With this expression for $r$, we obtain 
\begin{align*}
  \psi'(t) &= 1 - 2\jac{dn}{p}^2(\FAC\, t),
\end{align*}
from \eqref{eq:parPsiJEF} (up to sign). Integration provides the last piece in 
the proof of the following theorem. 

\begin{thm}[Classification of parabolic rotational sub-horospherical 
surfaces in $\H^3$]\label{thm:smallHpar}
  Every parabolic rotational constant mean curvature surface with
  $H^2<1$ in $\H^3 \subset \R^{3,1}$ is given by
    \begin{align*}
    (\theta, t)\mapsto \left(
      \frac{1+r^2(\psi^2 + \theta^2)}{2r},
      r(t)\theta, 
      r(t), 
      r(t)\psi(t)
    \right),
  \end{align*}
  and the functions
  \begin{align*}
    r(t) &= \frac{\sqrt{C^+}}{\jac{cn}{p}\left(
    \FAC\,t
    \right)}, \textrm{ with } p^2 = \left\{\begin{array}{ll}
      \frac{1+H}{2} &\textrm{for } C>0 \\
      \frac{1-H}{2} &\textrm{for } C<0.
    \end{array}\right.\\
    \psi(t) &=
      s - \frac{2}{\FAC} \left(E_p \circ \jac{am}{p}\right)(\FAC\, t),
  \end{align*}
  where $\FAC = \sqrt{|2C|}$.
\end{thm}

\begin{figure}%
\centering
\begin{tikzpicture}
	\node at (-5.5,0) {\includegraphics[width=0.3\textwidth]{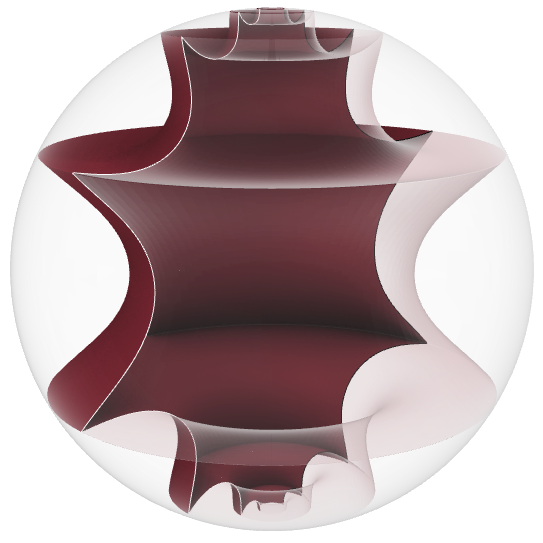}};
	\node at (0,0) {\includegraphics[width=0.3\textwidth]{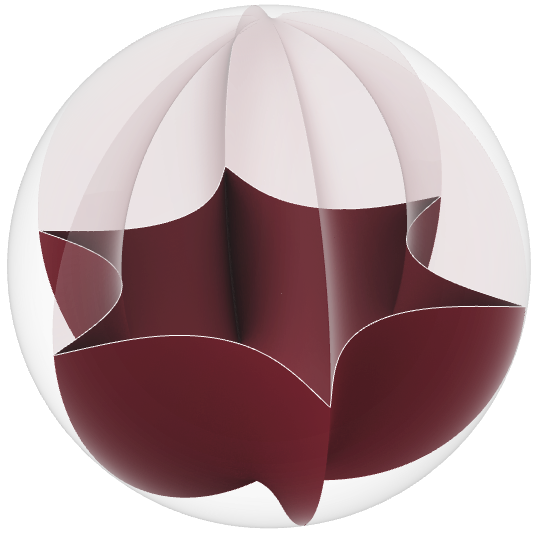}};
	\node at (5.5,0) {\includegraphics[width=0.3\textwidth]{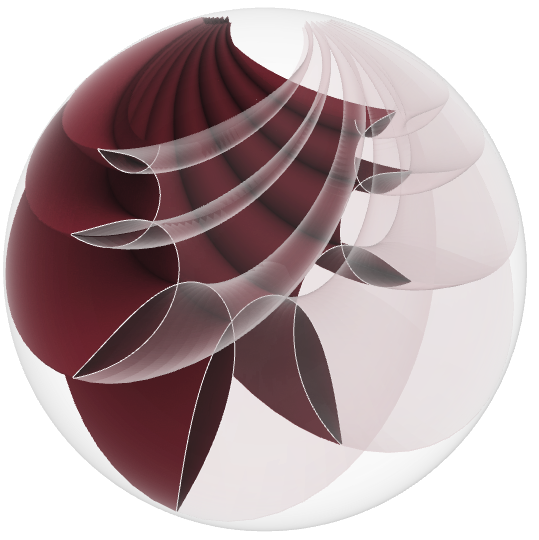}};
\end{tikzpicture}
\caption{Surfaces with cmc $H=0.3$ obtained using the parametrizations of 
\Cref{thm:smallH,thm:smallHpar} (elliptic, hyperbolic and parabolic 
rotation from left to right). Choosing the right parameter for the 
hyperbolic rotational case, such that the profile curve closes, amounts
to a problem similar to the one discussed in \Cref{thm:cmcTori}. Plotted 
and rendered using Rhino and mathematica.}%
\label{}%
\end{figure}

%% file: Content/ChannelLW.tex
In \cite{hertrich-jeromin2023} a complete and transparent classification of 
channel linear Weingarten (clW) surfaces in space forms in terms of Jacobi
elliptic functions was given for 
\begin{itemize}
  \item all clW surfaces in $\S^3$ (\cite[Thm 6.2]{hertrich-jeromin2023}),
  \item ``most'' clW surfaces in $\H^3$ (%
    \cite[Thm 6.10]{hertrich-jeromin2023}).
\end{itemize}
The surfaces missing in $\H^3$ were those in parallel families that do not 
contain surfaces of constant Gauss curvature. We will close this gap 
here. We will state a Bonnet-type theorem that classifies lW surfaces
in $\H^3$ in terms of constant curvature surfaces in their parallel families 
(\Cref{thm:Bonnet}) and then collect the needed parametrizations to complete
\cite[Thm 6.10]{hertrich-jeromin2023}. For simplicity, we restrict 
our attention to $\kappa=-1$, but the results extend to other 
hyperbolic space forms.

Recall that we call a surface lW, if its Gauss and mean 
curvature satisfy a (non-trivial) relationship of the form
\begin{align}\label{eq:LW2}
  aK + 2bH + c = 0,
\end{align}
for constants $a, b, c \in \R^3$. Because the principal curvatures can be 
expressed in terms of the curvature spheres $s_i$ via 
$\kappa_i = \tfrac{(\s_i, \q)}{(\s_i, \p)}$, this can be written as 
\begin{align*}
  \left((\s_1,\q), (\s_1, \p)\right) W 
  \left((\s_2,\q), (\s_2, \p)\right)^\tra
  =0, \textrm{ where } W=\begin{pmatrix}
    a &b \cr b &c
  \end{pmatrix}.
\end{align*}
As was described in \cite[Sect 4.8]{burstall2018}, parallel transformations of 
surfaces in a space form can be viewed as application of rotations $P^t$ in the 
plane $\spn{\p, \q}$, resulting in a change of the linear Weingarten condition
of the form
\begin{align*}
  W \mapsto P^tW(P^t)^\tra. 
\end{align*}
In the specific case of $\H^3$, we have
\begin{align*}
  P^t=\begin{pmatrix}
    \cosh(t) &\sinh(t) \cr \sinh(t) &\cosh(t)
  \end{pmatrix}.
\end{align*}
This can be used to prove the following proposition (for a proof see 
\cite[Sect~4.8]{burstall2018} or \cite[Thm~3.15]{polly2022}).

\begin{prop}\label{thm:Bonnet}
  Let $f:\Sigma^2 \to \H^3$ be a non-tubular linear Weingarten surface 
  satisfying \eqref{eq:LW2}. If $D = -\det W <0$, then $f$ is parallel
  to a surface of (negative) constant Gauss curvature. For $D>0$, there
  are three cases: 
  \begin{enumerate}
    \item If $\left|\frac{a+c}{2}\right|>|b|$, then the parallel family of $f$
      contains a surface of constant (positive) Gauss curvature. 
    \item If $\left|\frac{a+c}{2}\right|=|b|$, is of Bryant type and its 
      parallel family either consists of intrinsically flat surfaces ($K+1=0$)
      or contains either a surface of constant mean curvature $H=1$ or constant 
      harmonic mean curvature $\overline{H}=1$.
    \item If $\left|\frac{a+c}{2}\right|<|b|$, then the parallel family of $f$
      contains either a constant mean curvature $H<1$ surface or a constant
      harmonic mean curvature $\overline{H}>1$ surface. 
  \end{enumerate}
\end{prop}

This proposition (together with \Cref{thm:CLWIsRot}) can be used to classify
channel linear Weingarten surfaces with explicit parametrizations up to parallel
transformation. Relying on the results of \cite{hertrich-jeromin2023} and  
Section~\ref{sec:channelCMC}, we still have to give explicit parametrizations of
constant harmonic mean curvature (chc) $|\overline{H}|\geq 1$ surfaces. 

For a chc surface $(\Lambda, \p^\pm)$ the sphere complexes $\p^\pm$ take the 
form
\begin{align}\label{eq:HMCQuantities}
  \p^+ = \q ,\quad \p^- = \p-\frac{1}{\overline{H}}\q,
\end{align}
and Christoffel dual lifts of the isothermic sphere congruences $s^\pm$ are 
given as combinations of the space form lifts $\f, \n$ of the point sphere 
envelope and the tangent plane congruence via
\begin{align*}
  \s^+ = \n, \quad \s^- = \f + \frac{1}{\overline{H}} \n,
\end{align*}
see \Cref{bsp:conservedQuantities}. Note that $s^+$ is the tangent plane 
congruence of $\Lambda$ (which corresponds to the hyperbolic Gauss map of the 
surface, see \cite{cecil2008}) and $s^-$ is the middle sphere congruence (see 
\cite[\textsection 67]{blaschke1929}). One could use the equations developed in 
Subsection~\ref{sec:rotLW} to parametrize $s^+$ and $s^-$ and look for the 
space form lift. However, noting that $\f$ and $\n$ satisfy
\begin{align*}
  (\f,\n) = (\f,\n_t) = (\f,\n_\theta) = 0,
\end{align*}
we can determine $\f$ as the re-scaled ($\R^{3,1}$-)cross product of 
$\n, \n_t$ and $\n_\theta$:
\begin{align*}
  (u\times v\times w) = \det
  \begin{pmatrix}
    -\e_0 &\e_1 &\e_2 &\e_3 \\
    u_0 &u_1 &u_2 &u_3 \\
    v_0 &v_1 &v_2 &v_3 \\
    w_0 &w_1 &w_2 &w_3
  \end{pmatrix},
\end{align*}
for an orthonormal basis $(\e_0, \e_1, \e_2, \e_3)$ of $\R^{3,1}$ with timelike
$\e_0$. 

From \eqref{eq:HMCQuantities} we obtain
\begin{align*}
  \kappa^+ = 1,~\kappa^- = \frac{1}{\overline{H}^2} - 1,~\varepsilon = -\frac{1}{\overline{H}}
  ~\Rightarrow~\Delta = -1.
\end{align*}
As in the previous section, we equip $\R^{3,1}$ with an orthonormal basis 
$(\e_0, \e_1, \e_2, \e_3)$ for non-parabolic rotations and with a 
pseudo-orthonormal basis $(\v, \e_1, \o, \e_2)$ for parabolic rotations. We use 
this in \eqref{eq:rJEF} and \eqref{eq:psiJEF} (and \eqref{eq:parRJEF}, 
\eqref{eq:parPsiJEF}) to obtain the following theorem.

\begin{bsp}\label{bsp:BryantCHC}
  As in Section~\ref{sec:channelCMC}, we treat the case $\overline{H}=1$--where
  the polynomial $p_C$ degenerates--separately. In this
  instance, we learn (for non-parabolic rotations)
  \begin{align*}
    r(s) = \sqrt{\frac{C^2}{A}} \sinh(\sqrt{A}~s),
  \end{align*}
  where $A=\tfrac{1}{\kappa_1} +2C$ ($A=2C$). Note that, 
  contrary to the cmc case, negative values of $A$ also yield a real solution. 
\end{bsp} 

\begin{thm}[Classification of constant harmonic mean curvature 
$\overline{H}^2 > 1$ surfaces in $\H^3$]\label{thm:chc}
  The Gauss map of every rotational constant harmonic mean curvature surface 
  with $\overline{H}^2 > 1$ in $\H^3 \subset \R^{3,1}$ can be parametrized 
  by
  \begin{itemize}
    \item[elliptic rotations: $\kappa_1 = 1$] \begin{align*}
        (\theta, t)\mapsto \left(
          \sqrt{|r(t)^2-1|} \cosh \psi(t), 
          \sqrt{|r(t)^2-1|} \sinh \psi(t),
          r(t) \cos\theta, 
          r(t) \sin\theta
        \right),
    \end{align*}
    \item[hyperbolic rotations: $\kappa_1 = -1$] \begin{align*}
        (\theta, t) \mapsto \left(
          r(t) \cosh\theta, 
          r(t) \sinh\theta, 
          \sqrt{r^2(t)+1} \cos\psi(t), 
          \sqrt{r^2(t)+1} \sin\psi(t)
        \right)
    \end{align*}
      with
      \begin{align*}
        r(t) &= \sqrt{C^+}\jac{cn}{p}(\FAC\,t),~ p^2=\frac{C^+}{C^+ - C^-}\\
        \psi(t) &=\frac{1}{\overline{H}}
        \left\{
          -\kappa_1 t + \frac{\overline{H}C - \kappa_1}{\kappa_1C^+ - 1}
          \frac{1}{\FAC}\Pii{p}{\frac{\kappa_1C^+}{\kappa_1 C^+ - 1}}{\FAC\,t}
        \right\},
      \end{align*}
      where
      $\FAC=\sqrt{1-\tfrac{1}{\overline{H}^2}}\sqrt{C^+ - C^-}$ and
      $C^+ > 0 > C^-$ are the (real) roots of
      \begin{align}\label{eq:chcPoly}
        p_C(x) = \left(\frac{1}{\overline{H}^2}-1\right)x^2 + 
          \left(\kappa_1 - \frac{2C}{\overline{H}}\right)x + C^2.
      \end{align}
    \item[parabolic rotations:] \begin{align*}
      (\theta, t)\mapsto \left(
      \frac{-1+r^2(\psi^2 + \theta^2)}{2r},
      r(t)\theta, 
      r(t), 
      r(t)\psi(t)
    \right),
    \end{align*}
    with
    \begin{align*}
      r(t) &=\sqrt{C^+}\jac{cn}{p}(\FAC\,t),~ p^2=\frac{C^+}{C^+ - C^-} \\
      \psi(t) &= -\frac{1}{\overline{H}}t +
        \frac{C}{C^+} \frac{1}{\FAC}
        \Pii{p}{1}{\FAC\,t},
    \end{align*}
    where $\FAC=\sqrt{2|C|}$ and
    \begin{align*}
      C^\pm = \frac{C}{\tfrac{1}{\overline{H}} \pm \operatorname{sign}C}.
    \end{align*}  
  \end{itemize}
\end{thm}

\begin{rem}
  In the case of elliptic rotations, the polar coordinates used in $\Pi_2$ 
  are not necessarily causal (see Remark \ref{rem:polarCoordinates}). However,
  this does not effect the solution of \eqref{eq:rJEF}, only how we compute 
  $d$. This explains the absolute values in the parametrization given for 
  this case. 
\end{rem}

\begin{proof}
Using \eqref{eq:rJEF} and \eqref{eq:psiJEF} (or \eqref{eq:parRJEF} and 
\eqref{eq:parPsiJEF} for the parabolic case), we obtain a parametrization 
of $s^+$, i.e., the tangent plane congruence. The computation is very 
similar to the cmc $H<1$ case, described in Subsection~\ref{sec:smallH}. 
We therefor only give a basic outline of this step and refer for details 
to the proof of \Cref{prop:rOfSmallH}. Given the specific parametrizations,
we assume $\kappa_1 = -\kappa_2 \in \{\pm 1\}$ for the non-parabolic case.

First, the polynomial \eqref{eq:polynomJEF} takes the form \eqref{eq:chcPoly}.
Its roots, given by
  \begin{align*}
    C^\pm = \frac{
      \kappa_1 - \frac{2C}{\overline{H}} \mp \sqrt{
        1-4\frac{\kappa_1 C}{\overline{H}} + 4C^2
      }
    }{2\left(\frac{1}{\overline{H}^2-1}\right)},
  \end{align*}
are real for all values of $C$ and satisfy $C^+\geq 0\geq C^-$. This 
holds \emph{mutatis mutandis} for the polynomial \eqref{eq:parPolynomJEF} in the
parabolic case. 
  
The solution of \eqref{eq:rJEF} (or \eqref{eq:parRJEF}) is given by\footnote{
  Note that $\FAC^2$ is non-negative because $\kappa- < 0$, which is the key 
  difference to the sub-horospherical case investigated in 
  Subsection~\ref{sec:smallH}.
  }
  \begin{align*}
    r(t) = \sqrt{C^+}\jac{cn}{p}\left(\FAC\,t
    \right),
    \textrm{ with } p^2=\frac{C^+}{C^+ - C^-}.
  \end{align*}
This solution and \eqref{eq:psiJEF}, \eqref{eq:parPsiJEF} yield
  \begin{align*}
  \begin{array}{ll}
    \textit{non-parabolic:} &\psi'(t) =\frac{1}{\overline{H}}
      \left\{
        -\kappa_1 + \frac{\overline{H}C-\kappa_1}{\kappa_1 C^+-1}
        \frac{1}{1-\frac{\kappa_1 C^+}{\kappa_1 C^+ -1}\jac{sn}{p}^2(\FAC\,t)}
      \right\}, \\
    \textit{parabolic:} &\psi'(t) = -\frac{1}{\overline{H}} + 
    \frac{C}{C^+}\frac{1}{1-\jac{sn}{p}^2(\FAC\,t)}.
  \end{array}
  \end{align*} 
  which integrates to the desired solution. 
\end{proof}

\begin{rem}
  Due to the duality mentioned in \cite[Sec 2]{burstall2018},
  the Gauss map $\n$ of a chc surface $\f$ is itself a spacelike cmc surface
  in de Sitter space. An investigation of spacelike lW surfaces 
  in Lorentz space forms shall be carried out in a future project.
\end{rem}

Together with \cite[Thm 6.10]{hertrich-jeromin2023}, and Examples \ref{bsp:BryantCMC} 
and \ref{bsp:BryantCHC} the following theorem 
provides a complete classification of channel linear Weingarten surfaces in
hyperbolic space $\H^3$. 

\begin{thm}\label{thm:cLWH3}
  Every channel linear Weingarten surface in $\H^3$ with $|\frac{a+c}{2}|<|b|$ 
  is parallel to a rotational surfaces determined by one of the 
  parametrizations given in \Cref{thm:smallH,thm:smallHpar,thm:chc}. 
\end{thm}